\documentclass[10pt]{amsart}
\usepackage{mathtools}
\usepackage{kpfonts}
\usepackage{amsmath}
\usepackage{amsthm}
\usepackage{amssymb}
\usepackage{mathdots}
\usepackage{pst-node}
\usepackage{enumitem}
\usepackage{algpseudocode}
\usepackage{algorithmicx}

\newtheorem{theorem}{Theorem}[section]
\newtheorem{lemma}[theorem]{Lemma}
\newtheorem{proposition}[theorem]{Proposition}

\newtheorem{example}[theorem]{Example}

\theoremstyle{remark}
\newtheorem{remark}[theorem]{Remark}

\numberwithin{equation}{section}

\def\dim{\mathop{\rm dim}\nolimits}

\def\Ima{\mathop{\rm Im}\nolimits}
\def\Orb{\mathop{\rm Orb}\nolimits}

\def\Her{\mathop{\rm Her}\nolimits}

\parskip=\smallskipamount

\makeatletter
\newcommand{\doublewidetilde}[1]{{%
  \mathpalette\double@widetilde{#1}%
}}
\newcommand{\double@widetilde}[2]{%
  \sbox\z@{$\m@th#1\widetilde{#2}$}%
  \ht\z@=.9\ht\z@
  \widetilde{\box\z@}%
}
\makeatother

\makeatletter
\renewcommand*\env@matrix[1][*\c@MaxMatrixCols c]{%
  \hskip -\arraycolsep
  \let\@ifnextchar\new@ifnextchar
  \array{#1}}
\makeatother

\begin{document}

\title{Isotropy groups of the action of orthogonal *congruence on Hermitian matrices}
\author{Tadej Star\v{c}i\v{c}}
\address{Faculty of Education, University of Ljubljana, Kardeljeva Plo\v{s}\v{c}ad 16, 1000 Lju\-blja\-na, Slovenia}
\address{Institute of Mathematics, Physics and Mechanics, Jadranska
  19, 1000 Ljubljana, Slovenia}
\email{tadej.starcic@pef.uni-lj.si}
\subjclass[2000]{15A24, 15B57, 51H30}
\date{May 23, 2023}


\keywords{
isotropy groups, complex orthogonal matrix, Hermitian matrix, matrix equation, Toeplitz matrix\\
\indent Research supported by grants P1-0291 and J1-3005
from ARRS, Republic of Slovenia.}

\begin{abstract} 
We present a procedure which enables
the computation and the description of structures of isotropy subgroups of the group of complex orthogonal matrices with respect to the action of *congruence on Hermitian matrices. 
A key ingredient in our proof is an algorithm giving solutions of a certain rectangular block (complex-alter\-na\-ting) upper triangular Toeplitz matrix equation. 
\end{abstract}

\maketitle

\section{Introduction}

We denote by $\Her_n$ the real vector space of all $n$-by-$n$ Hermitian matrices; $A$ is Hermitian if and only if $A=A^{*}:=(\overline{A})^{T}$. Let further $O_n(\mathbb{C})$ be the group of complex orthogonal
$n$-by-$n$ 
matrices.
A matrix $Q$ is orthogonal precisely when $Q^{-1}=Q^{T}$. 
The action of \emph{orthogonal *congruence}
on $\Her_n$ is defined as follows:
\begin{align}\label{aos}
\Phi\colon O_n(\mathbb{C})\times \Her_n\to \Her_n,\qquad (Q,A)\mapsto Q^{*}AQ.
\end{align}
%
The study of Hermitian matrices under *congruence is indeed quite general, as can be concluded from Hua's fundamental result \cite[Theorem 12]{Hua45} on the geometry of Hermitian matrices (extended by Wan \cite[Theorem 6.4]{Wan}); see also the paper by Radjavi and \v {S}emrl \cite{RadjaviSe}. 
On the other hand
(\ref{aos}) can be seen as a representation of $O_n(\mathbb{C})$ as a real classical group (e.g. see the monograph \cite{Weyl}).

The \emph{isotropy group} at $A\in \Her_n$ 
with respect to the action (\ref{aos}) is denoted by
\begin{equation}\label{isog}
\Sigma_A:=\{Q\in O_n(\mathbb{C})\mid Q^{*}AQ=A\}.
\end{equation}
Isotropy groups provide an important information about a group action (see textbooks \cite{GOV,Milne}).
In a generic case (on a complement of a real analytic subset of codimension $1$) 
isotropy groups for (\ref{aos}) are clearly trivial (Proposition \ref{stabs11}), while 
in general the situation is more involved. The main purpose of this paper is to give an inductive procedure that enables the computation and the description of a group structure of an isotropy group (\ref{isog}) in a nongeneric case (Theorem \ref{stabw} and Theorem \ref{stabz}). Analoguous results for skew-Hermitian matrices under orthogonal *conjugation are valid as well.
Key ingredients in the proof 
are Lemma \ref{EqT} and Lemma \ref{EqTca}.
They provide solutions of certain block rectangular (complex-alternating) upper triangular Toep\-litz matrix equations.
These equations characterize orthogonality of a solution $Q$ of the
equation $A\overline{Q}=QA$ (or equivalently $A=QAQ^{*}$, 
i.e. $Q^{*}\in \Sigma_A$); for a general $Q$ this equation was 
solved by Bevis, Hall and Hartwig \cite{BHH}.

In contrast to the complex case, the situation in the real case is simple. 
Each real symmetric matrix 
is real orthogonally similar to $\Lambda=\oplus_{r=1}^{N} (\oplus_{j=1}^{m_r} \lambda_j )$ with $\lambda_1,\ldots,\lambda_N\in \mathbb{R}$ pairwise distinct.
Since 
$Q^{T}\Lambda Q=\Lambda$ for real orthogonal $Q$ transforms to the Sylvester equation $\Lambda Q=Q \Lambda$, 
the isotropy group at $\Lambda$ with respect to real orthogonal similarity consists of matrices $Q=\oplus_{r=1}^{N} Q_r$ with $Q_r$ real orthogonal of size $m_r\times m_r$.

Pairs $(A,B)$ with $A$ arbitrary and $B$ symmetric (i.e. $B=B^{T}$) 
with respect to transformations $(cP^{*}AP, P^{T}BP)$ for a nonsingular matrix $P$ and $c\in \mathbb{C}\setminus\{0\}$
are studied in CR-geometry in the theory of CR-singularities of codimension $2$.
Normal forms under this action for $2\times 2$ matrices 
were obtained
by Coffman \cite{Coff}.
In higher dimensions
the isotropy groups of (\ref{aos}) are expected to some extend to be applied to tackle this problem as well as a closely related problem of simultaneous reduction of $(A,B)$ under transformations $(P^{*}AP, P^{T}BP)$ with $P$ nonsingular.
By applying Autonne-Takagi factorization we first reduce $(A,B)$ to $(A',I)$ with the identity $I$. Next, we write $A'=H_1+iH_2$ with $H_1,H_2$ Hermitian.
We put $H_1$ 
into Hong's orthogonal *congruence normal form \cite{Hong89}
and then simplify $H_2$ by
using matrices from the isotropy group $\Sigma_{H_1}$, as they keep $H_1$, $I$ intact. 
We add that 
a reduction of Hermitian-symmetric pairs
was considered by Hua \cite{Hua46}, Hong \cite{Hong89}, Hong, Horn and Johnson \cite{HongHornJohn},
among others.

%
%
%
%
\section{The main results}\label{secIG}

Isotropy groups corresponding to elements of $\Orb (A):=\{Q^{*}AQ \mid Q\in O_n(\mathbb{C})\}$ (i.e. the orbit of $A$ with respect to (\ref{aos})) are conjugate,
thus 
it suffices to com\-pute them
for representatives of orbits. 
Hong \cite[Theorem 2.7]{Hong89} pro\-ved that each Hermitian matrix $A$
is orthogonally *congruent to a matrix of the form
%
\begin{equation}\label{NFE}
\mathcal{H}^{\varepsilon}(A)=\bigoplus_{j}^{}\varepsilon_j H_{\alpha_j}(\lambda_j)\oplus\bigoplus_{k}^{} K_{\beta_k}(\mu_k)\oplus\bigoplus_{l}^{} L_{\gamma_l}(\xi_l),
\end{equation}
in which $ \lambda_j\geq 0$, $\mu_k>0$, $\Ima (\xi_l)>0$, 
$\varepsilon=(\varepsilon_1,\varepsilon_2,\ldots)$, all $\varepsilon_j\in\{1,-1\}$ with $\varepsilon_j=1$ if
$\lambda_j=0$ and $\alpha_j$ odd, and 
where $\lambda_j^{2}$, $-\mu^{2}_k$ and $\xi_l^{2}$ are nonnegative, positive and nonreal eigenvalues of $A\overline{A}$, respectively;
%
\vspace{-1mm}
\small
\begin{align}\label{Hmz}
&H_n(z):=
\frac{1}{2}\left(
\begin{bmatrix}
0  &                &    1  & 2z\\
 &   \iddots           &   \iddots     & 1 \\
1   &  \iddots  & \iddots &  \\
2z   & 1          &     & 0 \\
\end{bmatrix}
+
i
\begin{bmatrix}
0  & 1 &             & 0 \\
-1 & \ddots &   \ddots        &   \\
   & \ddots  &  \ddots & 1 \\
0   &          & -1    & 0 \\
\end{bmatrix}
\right)\qquad (n\textrm{-by-}n),\\
%
& \label{KLmz}
K_{n}(z):=
\begin{bmatrix}
0       & -iH_n(z)    \\
iH_n(z) &  0    \\
\end{bmatrix}, \qquad
L_{n}(z):=
\begin{bmatrix}
0       & H_n(z)    \\
H^{*}_n(z) &  0    \\
\end{bmatrix}.
\end{align}
\normalsize
The author \cite[Theorem 1.1]{TSH} showed that (\ref{NFE}) is 
uniquely determined up to a permutation of its blocks.
We add that *congrunce canonical forms and dimensions of their orbits for general matrices are known as well (\cite{HornSergei}, \cite{TeranDopi1}), and that 
the isotropy subgroups of invertible integer matrices under congruence 
at
symmet\-ric Gram matrices  
of edge-bipartite graphs 
were 
studied in
\cite{MM},  
\cite{S3}.

By applying 
results from \cite[Sec. 2]{BHH} on solutions of the equation $A\overline{Y}=YA$,
we immediately conclude the following facts; check also Proposition \ref{resAoXXA} (\ref{resAoXXA1}).

\begin{proposition} \label{stabs11}
\begin{enumerate}[label= \arabic*.,ref=\arabic*.,
leftmargin=20pt]
\item 
Let $\rho_1,\ldots,\rho_n\in \mathbb{C}$ be all distinct and let 
$
\mathcal{H}^{\varepsilon}=\bigoplus_{j=1}^{n}\mathcal{H}_j^{\varepsilon}$ be of the form (\ref{NFE}),
in which $\mathcal{H}_j^{\varepsilon}$ is a direct sum whose summands 
correspond to the eigenvalue $\rho_j$ of $\mathcal{H}^{\varepsilon}\overline{\mathcal{H}^{\varepsilon}}$. 
Then $\Sigma_{\mathcal{H}^{\varepsilon}}=\bigoplus_{j=1}^{n}\Sigma_{\mathcal{H}^{\varepsilon}_j}$.
\vspace{-2mm}
\item
If 
$\mathcal{H}^{\varepsilon}= 
\bigoplus_{j=1}^{n}\varepsilon_j \lambda_{j} \oplus \bigoplus_{l=1}^{m} 
\begin{bsmallmatrix}
0 & \xi_l \\
\overline{\xi}_l & 0
\end{bsmallmatrix}
$ (a generic canonical form), in which $\lambda_{j}\geq 0$, $\xi_l\in \mathbb{C}\setminus{\mathbb{R}}$ 
are all distinct constants and all $\varepsilon_j\in\{1,-1\}$, 
then $\Sigma_{\mathcal{H}^{\varepsilon}}$ is trivial.
\end{enumerate}
\end{proposition}

In Sec. \ref{notation},
we describe
nonsingular solutions of 
$A\overline{Y}=YA$
by the follow\-ing matrices. 
Given $\alpha=(\alpha_1,\ldots,\alpha_N)$ with $\alpha_{1}>\ldots >\alpha_{N}$ and $\mu=(m_1,\ldots,m_N)$ let $\mathbb{T}^{\alpha,\mu}$ and $\mathbb{T}_c^{\alpha,\mu}$ consist of
$N$-by-$N$ block matrices with $\alpha_r$-by-$\alpha_s$ blocks of the form:
\vspace{-1mm}
\begin{align}\label{0T0}
\mathcal{X}=[\mathcal{X}_{rs}]_{r,s=1}^{N},\qquad
&\mathcal{X}_{rs}=
\left\{
\begin{array}{ll}
[0\quad \mathcal{T}_{rs}], & \alpha_r<\alpha_s\\
\begin{bmatrix}
\mathcal{T}_{rs}\\
0
\end{bmatrix}, & \alpha_r>\alpha_s\\
\mathcal{T}_{rs},& \alpha_r=\alpha_s
\end{array}\right., \qquad 
b_{rs}:=\min\{\alpha_s,\alpha_r\},
\end{align}
in which $\mathcal{T}_{rs}=T(A_0^{rs},\ldots,A_{b_{rs}-1}^{rs})$ and $\mathcal{T}_{rs}=T_c(A_0^{rs},\ldots,A_{b_{rs}-1}^{rs})$, respectively; $A_n^{rs}\in \mathbb{C}^{m_r\times m_s}$ and all $A_0^{rr} $ are nonsingular.
We use the standard notation 
$ \mathbb{C}^{m\times n}$ to denote the set of $m$-by-$n$ matrices, and 
let a $\beta$-by-$\beta$ \emph{block upper triangular Toeplitz}  
and a $\beta$-by-$\beta$ block \emph{complex-alternating upper triangular Toeplitz} matrix be:
\vspace{-1mm}
\setlength{\arraycolsep}{4.4pt}
\small
\begin{align*}
T(A_0,\ldots,A_{\beta-1}):=\begin{bmatrix}
  A_{0} & A_{1}                       & \ldots &    A_{\beta-1}  \\
0       & \ddots   &  \ddots    & \vdots \\
 \vdots &   \ddots           & \ddots   &   A_1 \\
0              & \ldots  &  0      & A_0
\end{bmatrix},
T_c(A_0,\ldots,A_{\beta-1}):=\begin{bmatrix}
  A_{0} & A_{1}         &  \ldots             & \ldots &    A_{\beta-1}  \\
0       & \overline{A}_0 & \overline{A}_{1}   &      & \vdots \\
 \vdots & \ddots            & A_0              & \ddots &   \vdots \\ 
 \vdots &  & \ddots   & \ddots            &  \vdots\\
0       & \ldots            &  \ldots &  0      & \ddots
\end{bmatrix},
\end{align*}
\normalsize
%
respectively, in which $A_0 ,\ldots,A_{\beta-1}$ are of the same size, $T(A_0,\ldots,A_{\beta-1})=[T_{jk}]_{j,k=1}^{\beta}$,  $T_c(A_0,\ldots,A_{\beta-1})=[T_{jk}']_{j,k=1}^{\beta}$ with $T_{jk}=T_{jk}'=0$ for $j>k$ and with $T_{(j+1)(k+1)}=T_{jk}$, $T_{(j+1)(k+1)}'=\overline{T}_{jk}'$. When in addition $A_0$ is the 
identity matrix,
they are called block (complex-alternating) upper \emph{unitriangular} Toeplitz. 

\begin{example}\label{exX}
Examples of matrices of the form (\ref{0T0}) are ($\alpha_1=3$, $\alpha_2=2$): 
\small
\[
\mathcal{X}=
\begin{bmatrix}[ccc|cc]
A_1 & B_1 & C_1 &  G_1  &  H_1 \\ 
0   & A_1 & B_1     &  0  &   G_1 \\
0   & 0   & A_1     &  0 & 0  \\ 
\hline
0     & N_1 & P_1                         &  A_2  &  B_2 \\  
0    & 0   & N_1              &  0    &   A_2    
\end{bmatrix} , \qquad
\mathcal{X}_c=
\begin{bmatrix}[ccc|cc]
A_1 & B_1 & C_1 &  G_1  &  H_1  \\  
0   & \overline{A}_1 & \overline{B}_1     &  0  &   \overline{G}_1 \\  
0   & 0   & A_1     &  0 & 0 \\  
\hline
0     & N_1 & P_1                         &  A_2  &  B_2 \\ 
0    & 0   & \overline{N}_1              &  0    &   \overline{A}_2  
\end{bmatrix}.
\]
\normalsize
\end{example}

Let $I_n$ be the $n$-by-$n$ identity matrix.
Given $g =I_p\oplus -I_q$ denote by $O_{p,q}(\mathbb{C})$ (by $O_{p,q}(\mathbb{R})$) and $U_{p,q}(\mathbb{C})$ the complex (real) pseudo-orthogonal and pseudo-unitay group, consisting of matrics of all complex (real) matrices $Q$ such that $Q^{-1}=gQ^{T}g$ and $Q^{-1}=gQ^{*}g$, respectively.

We state our first main result;
we prove it in Sec.\ref{sec2}.

\vspace{5mm}

\begin{theorem}\label{stabw}
For $\mu=(m_1,\ldots,m_N)$, $\alpha=(\alpha_1,\ldots,\alpha_N)$ 
with $\alpha_{1}>\ldots >\alpha_{N}$, and 
$\varepsilon=\{\varepsilon_{r,j}\}_{r=1,\ldots,N}^{j=1,\ldots,m_r}$ with all $\varepsilon_{r,j}\in \{1,-1\}$,
let
\vspace{-1mm}
\small
\begin{equation*}
\mathcal{H}^{\varepsilon}=
\bigoplus_{r=1}^{N}\big(\bigoplus_{j=1}^{m_r}M_j^{r}\big), \quad
M_j^{r}=
\left\{
\begin{array}{ll}
\hspace{-1mm}\varepsilon_{r,j}H_{\alpha_r}(\sqrt{\rho}),  &  \hspace{-1mm} \rho \geq 0\\
\hspace{-1mm}K_{\alpha_r}(\sqrt{-\rho}),   &  \hspace{-1mm} \rho<0\\
\hspace{-1mm}L_{\alpha_r}(\sqrt{\rho}),                  &   \hspace{-1mm} \rho\in \mathbb{C}\setminus \mathbb{R}
\end{array}
\right.\hspace{-2mm}, \quad
B_{r}:=
\left\{
\begin{array}{ll}
\hspace{-1mm}\oplus_{j=1}^{m_r}\varepsilon_{r,j},  & \hspace{-1mm}  \rho\geq 0\\
\hspace{-1mm}\rho I_{ m_r}\oplus I_{ m_r},         &  \hspace{-1mm}  \rho<0\\
\hspace{-1mm}I_{2m_r},                              &  \hspace{-1mm} \rho \in \mathbb{C}\setminus \mathbb{R}
\end{array}
\right.\hspace{-2mm},
\end{equation*}
\normalsize
i.e. 
$\mathcal{H}^{\varepsilon}\overline{\mathcal{H}^{\varepsilon}}$ has precisely one eigenvalue $\rho$.
Then the isotropy group 
$\Sigma_{\mathcal{H}^{\varepsilon}}$ is conjugate (hence isomorphic) to the subgroup 
$\mathbb{X}\subset\left\{ 
\begin{array}{ll}
\mathbb{T}^{\alpha,\mu}, & \rho> 0\\
\mathbb{T}^{\alpha,\mu}\oplus \overline{\mathbb{T}}^{\alpha,\mu}, & \rho\in \mathbb{C}\setminus \mathbb{R} \\
\mathbb{T}_c^{\alpha,\mu}, & \rho= 0\\
\mathbb{T}^{\alpha,2\mu}, & \rho< 0
\end{array}
\right.$, 
where $\mathbb{T}^{\alpha,\mu}$, $\mathbb{T}_c^{\alpha,\mu}$ and $\mathbb{T}^{\alpha,2\mu}$
are defined by (\ref{0T0}).
Furthermore, each 
$\mathcal{X}\in \mathbb{X}$ for $\rho\in \mathbb{R}$ and each $\mathcal{X}\oplus \overline{\mathcal{X}}\in \mathbb{X}$ for $\rho\in \mathbb{C}\setminus \mathbb{R}$, with $\mathcal{X}$ of the form (\ref{0T0}), satisfy the following properties:
%
\begin{enumerate}[label={(\alph*)},ref={\alph*},
itemindent=1pt,leftmargin=18pt]

\item \label{stabs0} The nonzero entries of $\mathcal{X}_{rs}$ ror $r,s\in \{1,\ldots,N\}$ with $r>s$ can be taken freely. 
If either $\rho \in \mathbb{R}\setminus\{ 0\}$
or $\rho= 0$ with $\alpha_r$ odd ($\rho =0$ with $\alpha_r$ even), then $(\mathcal{X}_{rr})_{11}=A_0^{rr}$
are pseudo-orthogonal with $(A_0^{rr})^{-1}=B_{r}(A_0^{rr})^{T}B_{r}$ (pseudo-unitary with $(A_0^{rr})^{-1}=B_{r}(A_0^{rr})^{*}B_{r}$), while for $\rho \in \mathbb{C}\setminus \mathbb{R}$, matrices $A_0^{rr}$ are orthogonal.

\item \label{stabs4} 
For $r\in \{1,\ldots,N\}$ with 
$\alpha_r\geq 2$, $j\in \{1,\ldots,\alpha_r-1\}$ we have 
$(\mathcal{X}_{rr})_{1(1+j)}=A_{j}^{rr}=A_0^{rr}B_{r}Z_j^{r}+D_j^{r}$ for arbitrarily chosen 
%
\small
$Z_j^{r}=
\left\{
\begin{array}{ll}
\hspace{-1mm}-(Z_j^{r})^{*}, & \hspace{-1mm} \scriptstyle{\alpha_r-j \textrm{ even}, \rho=0} \\
\hspace{-1mm}-(Z_j^{r})^{T}, & \hspace{-1mm} \textrm{otherwise}
\end{array}
\right.$
\normalsize
and for some $D_j^{r}
$ depending polynomially on $A_{j'}^{r'r'}$ with $j'\in \{0,\ldots,j-1\}$, $r'\in \{1,\ldots,r\}$ and on the nonzero entries of $\mathcal{X}_{rs}$ for $r>s$ (described in (\ref{stabs0})). 
%

\quad
The nonzero entries of $\mathcal{X}_{rs}$ for $r,s\in \{1,\ldots,N\}$ with $r<s$ are uniquely determined (polynomially)
by the entries of $\mathcal{X}_{rs}$ with $r\geq s$ (described 
above
).

\item If $\rho>0$ then all $\mathcal{X}_{rs}$ are real, while for $\rho<0$ the upper triangular parts $\mathcal{T}_{rs}=T(A_0^{rs},\ldots,A_{b_{rs}-1}^{rs})$ of $\mathcal{X}_{rs}$ consist of $2$-by-$2$ block matrices of the form:
%
\vspace{-2mm}
\begin{equation}\label{Trsvw0}
\hspace{-6mm}
A_n^{}=\begin{bsmallmatrix}
V_n^{rs} & W_n^{rs}\\
\rho\overline{W}_n^{rs}+\overline{W}_{n-1}^{rs} & \overline{V}_n^{rs}
\end{bsmallmatrix},\,\,
\text{
\small $V_n^{rs},W_n^{rs} \in \mathbb{C}^{m_r\times m_s}$, $n\in \{1,\ldots,b_{rs}-1\}$; $V_{-1}^{rs}=W_{-1}^{rs}=0$
}.
\end{equation}
%
\end{enumerate}

\vspace{-1mm}
\noindent
In particular,
\vspace{-1mm}
\small
\begin{align*}
\dim_{\mathbb{R}} (\Sigma_{\mathcal{H}_{}^{\varepsilon}})=\left\{
\begin{array}{ll}
\displaystyle\sum_{r=1}^{N} m_r\big(\tfrac{1}{2}\alpha_r(m_r-1)+\sum_{s=1}^{r-1}\alpha_s m_s\big),    &   \rho>0\\
\displaystyle \sum_{r=1}^{N}\big(\alpha_r m_r^{2} +2\sum_{s=1}^{r-1}\alpha_s m_r m_s\big) - \sum_{\alpha_r \textrm{ even}}\tfrac{\alpha_r}{2}m_r-\sum_{\alpha_r \textrm{ odd}}\tfrac{\alpha_r+1}{2} m_r, &   \rho=0\\
\displaystyle\sum_{r=1}^{N} m_r\big(\alpha_r (2m_r-1)+2\sum_{s=1}^{r-1}\alpha_s m_s\big),              &  \rho<0\\
\displaystyle\sum_{r=1}^{N}2 m_r\big(\tfrac{1}{2}\alpha_r(m_r-1)+\sum_{s=1}^{r-1}\alpha_s m_s\big),                     & \rho\in \mathbb{C}\setminus \mathbb{R}
\end{array}
\right..
\end{align*}
\normalsize
\end{theorem}

\begin{remark}\label{OpI}
An algorithm to compute matrices in Theorem \ref{stabw} (\ref{stabs4}) 
is pro\-vi\-ded as part of its proof, more precisely, by 
Lemma \ref{EqT} and Lemma \ref{EqTca}.
\end{remark}

The following significant examples 
of matrices  
satisfy Theorem \ref{stabw} (\ref{stabs0}), (\ref{stabs4}).

\begin{example} (\cite[Example 3.1]{TSOS})
Fix $B_{r}$ nonsingular symmetric and let $Z_{n}^{r}$ be any skew-symmetric matrix (i.e. $Z_{n}^{r}=-(Z_{n}^{r})^{T}$); all of size $m_r\times m_r$. 
We set $W_0^{r}:=0$ and
\small
\begin{align}\label{asZ}
&\mathcal{W}=\bigoplus_{r=1}^{N}T(I_{m_r},W_1^{r},\ldots,W_{\alpha_r-1}^{r}),\qquad
W_{n}^{r}:=\frac{1}{2}B_{r}^{-1}\big(Z_{n}^{r}-\sum_{j=1}^{n-1}(W_j^{r})^{T}B_{r}W_{n-j}^{r}
\big).
\end{align}
\end{example}
\normalsize

\begin{example}
For $r\in \{1,\ldots ,N\}$, $n\in \{1,\ldots ,\alpha_r-1\}$, we are given $B_{r}$ nonsingular real symmetric with $\mathcal{B}_n^{r}:=\oplus_{j=1}^{n}B_r$ and let $Z_{n}^{r}$ be any skew-symmetric matrix for $\alpha_r-n$ odd (skew-Hermitian for $\alpha_r-n$ even); all of size $m_r\times m_r$. Set:
\vspace{-1mm}
\small
\begin{align}\label{asZ2}
&\mathcal{W}=\bigoplus_{r=1}^{N}T_c(I_{m_r},W_1^{r},\ldots,W_{\alpha_r-1}^{r}), \,\,\,\,
W_{n}^{r}:=\frac{1}{2}B_{r}^{-1}
\Big( Z_{n}^{r}-
\hspace{-0.5mm}
\left\{\begin{array}{ll}
\hspace{-0.5mm}\mathcal{A}_{n-1}^{r}\mathcal{B}_{n-1}^{r}\mathcal{P}_{n-1}^{r}, & \hspace{-0.5mm} \alpha_r \textrm{ even}\\
\hspace{-0.5mm}\overline{\mathcal{A}}_{n-1}^{r}\mathcal{B}_{n-1}^{r}\mathcal{P}_{n-1}^{r} & \hspace{-0.5mm} \alpha_r \textrm{ odd}
\end{array}
\right.\hspace{-0.5mm}\Big) , \\
&\mathcal{A}_n^{r}:=
\left\{
\begin{array}{ll}
\begin{bsmallmatrix}
(W_1^{r})^{T},(\overline{W}_2^{r})^{T},\ldots,(\overline{W}_{n-1}^{r})^{T},(W_{n}^{r})^{T}
\end{bsmallmatrix}, & n \textrm{ odd}\\ 
\begin{bsmallmatrix}(W_1^{r})^{T},(\overline{W}_2^{r})^{T},\ldots,(W_{n-1}^{r})^{T},(\overline{W}_{n}^{r})^{T}
\end{bsmallmatrix}, & n \textrm{ even}
\end{array}
\right.,\quad
\mathcal{P}_{2n-1}^{r}:=
\begin{bsmallmatrix}
\overline{W}_{2n-1}^{r} \\ 
W_{2n-2}^{r} \\ 
\vdots  \\ 
W_{2}^{r}\\
\overline{W}_{1}^{r} 
\end{bsmallmatrix}, 
\mathcal{P}_{2n}^{r}:=
\begin{bsmallmatrix}
\overline{W}_{2n}^{r} \\ 
W_{2n-1}^{r} \\ 
\vdots  \\ 
\overline{W}_{2}^{r} \\ 
W_{1}^{r}
\end{bsmallmatrix};\nonumber
\end{align}
\normalsize
the entry in the $j$-th column of $\mathcal{A}_n^{r}$ is $(W_{j}^{kr})^{T}$ (and $(\overline{W}_{j}^{kr})^{T}$) for $j$ odd (even), 
and the entry in the $j$-th row of $\mathcal{P}_{n}^{ks}$ is $(W_{n-j+1}^{ks})^{T}$ (and $(\overline{W}_{n-j+1}^{ks})^{T}$) for $j$ even (odd), $\mathcal{P}_0^{r}:=0$.
\end{example}

\begin{example}
Let a matrix
of the form (\ref{0T0}) 
have the identity as principal subma\-trix, formed by all blocks except those at the $p$-th, the $t$-th columns and rows, i.e.
\vspace{-2mm}
\begin{align}\label{Hptk}
&\mathcal{T}_{rs}=\left\{
\begin{array}{ll}
\oplus_{j=1}^{\alpha_r}I_{m_r}, &  r=s,\\
0,                        &  r\neq s 
\end{array}
\right., \{r,s\}\not\subset\{p,t\}.
\end{align}
\vspace{-1mm}

In particular, given $B_{r}$ nonsingular symmetric of size $m_r$-by-$m_r$, $F \in \mathbb{C}^{m_p\times m_t}$ and $0\leq k \leq \alpha_t-1$, $r\in \{1,\ldots ,N\}$, in \cite[Example 3.2]{TSOS} we set:
%
\begin{align}\label{propS}
&\mathcal{T}_{rr}=T(I_{m_r},A_1^{rr},\ldots,A_{\alpha_r-1}^{rr}), \quad r\in \{p,t\},\qquad p< t, \nonumber\\
&
A_{j}^{pp}=\left\{\begin{array}{ll}
a_{n-1}B_{p}^{-1}(F^{T}B_{t}FB_{p}^{-1})^{n}B_0^{r}, & j=n(2k+\alpha-\beta)\\
0,                      & \textrm{otherwise}
\end{array}
\right.,
\nonumber \\
&
A_{j}^{tt}=\left\{\begin{array}{ll}
a_{n-1}B_{t}^{-1}(B_{t}FB_{p}^{-1}F^{T})^{n}B_{t}, & j=n(2k+\alpha-\beta)\\
0,                      & \textrm{otherwise}
\end{array}
\right., \\
& 
a_{n}:=-\frac{1}{2^{2n+1}(n+1)}
\binom{2n}{n}
,\qquad
\mathcal{T}_{tp}=
N_{\alpha_t}^{k}(F ), \qquad \mathcal{T}_{pt}=
N_{\alpha_t}^{k}\big(-B_{p}^{-1}F^{T}B_{p}\big),
\nonumber
\end{align}
%
%
where $N_{\beta}^{k}(X)$ is a $\beta\textrm{-by-} \beta$ block matrix with $X$ on the $k$-th upper diagonal
(the main diagonal for $k=0$) and zeros otherwise.
For example,
if $N=2$, $\alpha_1=4$, $\alpha_2=2$,  
$m_1=2$, $m_2=3$, 
$B_{1}=I_2$, $B_{2}=I_3$,
then $F\in \mathbb{C}^{2\times 3}$ and we obtain
\small
\begin{equation}\label{exXC}
\begin{bmatrix}[cccc|cc]
I_{2} & 0 & -\tfrac{1}{2}F^{T}F & 0  &  -F^{T}  &  0 \\
0   & I_{2} & 0   & -\tfrac{1}{2}F^{T}F  &  0  &   -F^{T}\\
0   & 0   & I_{2} & 0    &  0 & 0 \\
0   & 0   &  0  &  I_{2}   & 0 &0 \\
\hline
0   & 0   & F & 0                         &  I_3  & 0  \\
0   & 0   & 0   & F              &  0    &   I_3  \\
\end{bmatrix}.
\end{equation}
\normalsize
%

The other intrieguing choice, with 
$G:=\left\{
\begin{array}{ll}
F, & k+\alpha_t \textrm{ odd}\\
\overline{F}, & k+\alpha_t \textrm{ even}
\end{array}
\right.$, is
\small
\begin{align}\label{propS2}
&\mathcal{T}_{rr}=T_c(I_{m_r},A_1^{rr},\ldots,A_{\alpha_r-1}^{rr}), \quad r\in \{p,t\},\qquad p< t \nonumber\\
&\nonumber
A_{n(2k+\alpha_p-\alpha_t)}^{pp}=a_{n-1} B_{p}^{-1} \left\{\begin{array}{ll}
(G^{T}B_{t}GB_{p}^{-1})^{n}B_{p}, & \alpha_p,\alpha_t \textrm{ odd}\\
(G^{T}B_{t}\overline{G}B_{p}^{-1})^{n}B_{p}, & \alpha_p,\alpha_t \textrm{ even}\\
(G^{T}B_{t}GB_{p}^{-1})^{\overline{n}}B_{p}, & \alpha_p \textrm{ even},\alpha_t \textrm{ odd}\\
(G^{T}B_{t}\overline{G}B_{p}^{-1})^{\overline{n}}B_{p}, & \alpha_p \textrm{ odd},\alpha_t \textrm{ even}
\end{array}
\right.,\\
&
A_{n(2k+\alpha_p-\alpha_t)}^{tt}=a_{n-1}B_{t}^{-1} \left\{\begin{array}{ll}
(B_{t}FB_{p}^{-1}F^{T})^{n}B_{t}, & \alpha_p,\alpha_t \textrm{ odd}\\
(B_{t}FB_{p}^{-1}\overline{F}^{T})^{n}B_{t}, & \alpha_p,\alpha_t \textrm{ even}\\
(B_{t}FB_{p}^{-1}F^{T})^{\overline{n}}B_{t}, & \alpha_p \textrm{ even},\alpha_t \textrm{ odd}\\
(B_{t}FB_{p}^{-1}\overline{F}^{T})^{\overline{n}}B_{t}, & \alpha_p \textrm{ odd},\alpha_t \textrm{ even}
\end{array}
\right.,\\
&A_{j}^{tt}=0,\quad  A_{j}^{pp}=0,\qquad  j\neq n(2k+\alpha_p-\alpha_t),\quad 
a_{n}=-\frac{1}{2^{2n+1}(n+1)}
\binom{2n}{n}
\nonumber\\
&
\mathcal{T}_{tp}=
Nc_{\alpha_t}^{k}(F ),
\qquad
\mathcal{U}_{pt}=
Nc_{\alpha_t}^{k}\big(-B_{p}^{-1}G^{T}B_{t}\big),
\nonumber
\end{align}
\normalsize
in which 
\small
$X^{\overline{n}}:=\left\{
\begin{array}{ll}
X\overline{X}X\cdots \overline{X}X, & n \textrm{ odd}\\
X\overline{X}\cdots X\overline{X}, & n \textrm{ even}
\end{array}
\right.$
\normalsize
is the complex-alternating produkt of $n$ factors with $X$ as odd (with $\overline{X}$ as even) factor, 
and
$Nc_{\beta}^{k}(X)$ is a complex-alternating $\beta\textrm{-by-} \beta$ Toeplitz with $X, \overline{X},X,\ldots$ on the $k$-th upper diagonal
(the main diagonal for $k=0$) and zeros otherwise.
%
%
If $N=2$, $\alpha_1=4$, $\alpha_2=2$, 
$m_1=2$, $m_2=3$, 
$B_{1}=I_2$, $B_{2}=I_3$,  
we have (c.f. (\ref{exXC})):
\small
\[
\begin{bmatrix}[cccc|cc]
I_{2} & 0 & -\tfrac{1}{2}F^{*}F & 0  &  -F^{*}  &  0 \\
0   & I_{2} & 0   & -\tfrac{1}{2}F^{T}\overline{F}  &  0  &   -F^{T} \\
0   & 0   & I_{2} & 0    &  0 & 0 \\
0   & 0   &  0  &  I_{2}   & 0 &0 \\ 
\hline
0   & 0   & F & 0                         &  I_3  & 0 \\ 
0   & 0   & 0   & \overline{F}              &  0    &   I_3 
\end{bmatrix}
\]
\end{example}
\normalsize

We now exihibit the structure of isotropy groups; the proof is given in Sec.\ref{sec2}.

\begin{theorem}\label{stabz}
Let $\mathcal{H}^{\varepsilon}$, $\mathbb{T}^{\alpha,\mu}$, $\mathbb{T}^{\alpha,2\mu}$ and $\mathbb{T}_c^{\alpha,\mu}$ be as in Theorem \ref{stabw}. Then $\Sigma_{\mathcal{H}^{\varepsilon}}$ is isomorphic to a semidirect product:
\vspace{-1mm}
\[
\Sigma_{\mathcal{H}^{\varepsilon}}\cong\mathbb{O}\ltimes \mathbb{V},
\]
in which $\mathbb{O}$ and $\mathbb{V}$ are described as follows:
%
\begin{enumerate}[label={\bf (\Roman*)},ref={\Roman*},itemsep=1pt,itemindent=1pt,leftmargin=17pt]
%
\item \label{stabz01} 
Suppose 
$\mathcal{H}_{}^{\varepsilon}=
\bigoplus_{r=1}^{N}\left(\bigoplus_{j=1}^{p_r} H_{\alpha_r}(\lambda)\oplus \bigoplus_{j=1}^{q_r}- H_{\alpha_r}(\lambda)\right)$ for $\lambda \geq 0$, $m_r:=p_r+q_r$. 
%
\vspace{-2mm}
\begin{enumerate}[label={ \roman*.},ref={\roman*}, itemsep=0pt,itemindent=1pt,
leftmargin=20pt]
\item If $\lambda>0$, then
$\mathbb{O}\subset \mathbb{T}^{\alpha,\mu}$ consists of all matrices $\mathcal{Q}=\bigoplus_{r=1}^{N}\big(\bigoplus_{j=1}^{\alpha_r} Q_r\big)$ with $Q_{r}\in O_{p_r,q_r}(\mathbb{R})$, while $\mathbb{V}\subset \mathbb{T}^{\alpha,\mu}$ is generated by all real matrices of the form (\ref{asZ}) and of the form (\ref{0T0}) with (\ref{Hptk}), (\ref{propS}) for $B_{r}=I_{p_r}\oplus -I_{q_r}$. 

\item If $\lambda=0$ and for $\alpha_r$ odd $m_r=p_r$, then
$\mathbb{O}\subset \mathbb{T}_c^{\alpha,\mu}$ 
consists of all matrices $\mathcal{Q}=\bigoplus_{r=1}^{N} (Q_r\oplus \overline{Q}_r\oplus Q_r\oplus \cdots)$ 
with $Q_r\in O_{m_r}(\mathbb{C})$ for $\alpha_r$ odd and $Q_r\in U_{p_r,q_r}(\mathbb{C})$ for $\alpha_r$ even, while $\mathbb{V}\subset \mathbb{T}_c^{\alpha,\mu}$ is generated by matrices of the form (\ref{asZ2}) and of the form (\ref{0T0}) with (\ref{Hptk}), (\ref{propS2}) for $B_{r}=I_{p_r}\oplus -I_{q_r}$.
\end{enumerate}
\vspace{-2mm}
(The possible summands $\bigoplus_{j=1}^{0} \pm H_{\alpha_r}(\lambda)$ and $\pm I_0$ are left out.)

\item \label{stabz2} If $\mathcal{H}_{}^{\varepsilon}=\bigoplus_{r=1}^{N}\left(\bigoplus_{k=1}^{m_r} K_{\alpha_r}(\mu)\right)$ for $ \mu >0$,
then 
$\mathbb{O}\subset \mathbb{T}^{\alpha,2\mu}$ 
consists of all matrices $\mathcal{Q}=\bigoplus_{r=1}^{N}\big(\bigoplus_{j=1}^{\alpha_r} Q_r\big)$ 
such that $(i\mu I_{m_r}\oplus I_{m_r})Q_r(i\mu I_{m_r}\oplus I_{m_r})\in O_{m_j,m_j}(\mathbb{C})$, and such that 
$Q_r=
\begin{bsmallmatrix}
V_0^{rr} &  W_0^{rr} \\
-\mu^{2}\overline{W}_0^{rr} &  \overline{V}_0^{rr} 
\end{bsmallmatrix}$ for some $V_0^{rr},W_0^{rr}\in \mathbb{C}^{m_r\times m_r}$, 
while
each $\mathcal{V}\in \mathbb{V}\subset \mathbb{T}^{\alpha,2\mu}$ can be written as $\mathcal{V}=\mathcal{V}_0\prod_{j=1}^{n}\mathcal{V}_j$, where $\mathcal{V}_0=\bigoplus_{r=1}^{N}\mathcal{W}_r$ with $\mathcal{W}_r$ upper unitriangular Toeplitz  
and $\mathcal{V}_1,\ldots,\mathcal{V}_n$ of the form (\ref{Hptk}); all $\mathcal{V}_0,\mathcal{V}_1,\ldots,\mathcal{V}_n$ are of the form (\ref{0T0}) with (\ref{Trsvw0}) and satisfying (\ref{stabs0}), (\ref{stabs4})
for $\rho<0$ in Theorem \ref{stabw}.

\item \label{stabz3} If $\mathcal{H}_{}^{\varepsilon}=\bigoplus_{r=1}^{N}\left(\bigoplus_{l=1}^{m_r} L_{\alpha_r}(\xi)\right)$, $ \xi^{2}\in \mathbb{C}\setminus  \mathbb{R}$, then 
$\mathbb{O}\subset \mathbb{T}^{\alpha,\mu}$ 
consists of all matrices $\mathcal{Q}=\bigoplus_{r=1}^{N}\big(\bigoplus_{j=1}^{\alpha_r} Q_r\big)$ 
with $Q_r\in O_{m_j}(\mathbb{C})$, and $\mathbb{V}\subset \mathbb{T}^{\alpha,\mu}$ is generated by matrices of the form (\ref{asZ}) and of the form (\ref{0T0}) with (\ref{Hptk}), (\ref{propS}) for $B_{r}=I_{m_r}$.
\end{enumerate}
\vspace{-1mm}

\noindent
In particular, $\mathbb{V}$ is unipotent of order at most $\alpha_1-1$ (nilpotent of class $\leq \alpha_1$).
\end{theorem}

\begin{remark}
Isotropy groups for $A$ and $iA$
under orthogonal *conjugation 
coincide, thus analogues of Theorem \ref{stabz} and  Theorem \ref{stabz} for skew-Hermitian matrices are valid.
\end{remark}

%
%
%
%
\section{The matrix equation $A\overline{Y}=YA$}\label{notation}

Given a square matrix $A$ we consider the matrix equation
%
\begin{equation}\label{eqAoXXB}
A\overline{Y}=YA.
\end{equation}
For $Y=PXP^{-1}$ with $P$ nonsingular, (\ref{eqAoXXB}) transforms to $B\overline{X}=XB$ for $B=P^{-1}A\overline{P}$; such $A$ and $B$ are said to be {\it consimilar}.
Bevis, Hall and Hartwig \cite{BHH} used the canonical form under consimilarity, given by Hong and Horn \cite[Theorem 3.1]{HongHorn}, 
to reduce (\ref{eqAoXXB}) to Sylvester equations.
In a similar fashion we shall solve (\ref{eqAoXXB})
by using Hermitian Hong's consimilarity canonical form
(\ref{NFE}) for $\varepsilon=(1,1,\ldots)$ \cite[p. 3-4]{Hong90}.
Consimilarity canonical forms were first developed by Haantjes \cite{Haantjes}, Asano and Nakayama \cite{AsanoNakaya}, but these are not suitable to solve (\ref{eqAoXXB}).

Recall the classical result \cite[Ch. VIII]{Gant}
on 
solutions of a
Sylvester equation.

\begin{theorem}\label{lemas}
Given $\lambda_1, \lambda_2\in \mathbb{C}$, an $m$-by-$n$ matrix $Y$ satisfies the matrix equation 
\[
J_m(\lambda_1)X=XJ_n(\lambda_2),\quad \qquad  
J_{\alpha}(\lambda):=\begin{bsmallmatrix}
                                                      \lambda    &  1       & \;     & 0    \\
						      \;     & \lambda     & \ddots & \;    \\     
						      \;     & \;      & \ddots &  1     \\
                                                      0     & \;      & \;     & \lambda   
                                   \end{bsmallmatrix},\quad \lambda\in \mathbb{C}\quad (\alpha\textrm{-by-}\alpha),
\]
if and only if 
either $\lambda_1\neq \lambda_2$ and $X=0$,
or $\lambda_1=\lambda_2$ and 
\begin{equation}\label{QTY}
X=\left\{\begin{array}{ll}
\begin{bmatrix} 
0 & T
\end{bmatrix}, & m<n \\
\begin{bmatrix}
T\\
0
\end{bmatrix}, & m>n\\
T, & n=m
\end{array}
\right.,
\end{equation}
in which $T$ is an $\beta$-by-$\beta$ upper triangular Toeplitz matrix ($\beta=\min\{m,n\}$). 
\end{theorem}

\begin{lemma}\label{lemaBHH}
Given matrices $M$ and $N$, let us consider the following equation 
\vspace{-1mm}
\begin{equation}\label{eqAoX2}
M\overline{Y}=YN.
\end{equation}
Denote the $n$-by-$n$ backward identity matrix by $E_{n}$ (with ones on the anti-diagonal).
\begin{enumerate}
\vspace{-1mm}
\item \label{BHH1}
If $M$ and $N$ are of the form (\ref{Hmz}) or (\ref{KLmz}) and such that 
$M\overline{M}$ and $N\overline{N}$ correspond to different eigenvalues,
it then follows that $Y=0$.
%
\vspace{-1mm}
\item \label{BHH2} If 
$M=H_m(\lambda)$ and 
$N=H_n(\lambda)$ with 
$\lambda$ positive (zero), then 
$Y$ satisfies (\ref{eqAoX2}) if and only if
$
Y=P_{m}^{-1}XP_n
$,
in which $X$ is an $m$-by-$n$ matrix of the form (\ref{QTY}) for $T$ an $\beta$-by-$\beta$ real (complex-alternating) upper triangular Toeplitz matrix with $\beta=\min\{m,n\}$, and $P_{\alpha}:=\frac{1}{\sqrt{2}}e^{{\scriptscriptstyle -\frac{i\pi}{4}}}(I_{\alpha}+iE_{\alpha})$ for $\alpha\in \{m,n\}$. 
\vspace{-1mm}
\item \label{BHHK} If 
$M=K_m(\mu)$ and 
$N=K_n(\mu)$ with 
$\mu>0$, then 
$Y$ satisfies (\ref{eqAoX2}) if and only if
%
$Y=Q_m^{-1}V_m^{-1}S_m(\mu)XS_n^{-1}(\mu)V_nQ_n$, 
in which $Q_{\alpha}:=e^{{\scriptscriptstyle \frac{i\pi}{4}}}(P_{\alpha}\oplus P_{\alpha})
$, \hspace{-0.5mm}
$V_{\alpha}:=e^{i\frac{\pi}{4}}(W_{\alpha}\oplus \overline{W}_{\alpha})$ with $W_{\alpha}:= \oplus_{j=0}^{\alpha-1} i^{j}$ for $\alpha\in \{m,n\}$, 
and
\begin{equation}\label{QTC}
X=\begin{bmatrix}
X_1 & X_2\\
J_m(-\mu^{2})\overline{X}_2 & \overline{X}_1
\end{bmatrix},
\end{equation}
where $X_1$, $X_2$ are $m$-by-$n$ matrices of the form (\ref{QTY}) for an $\beta$-by-$\beta$ upper triangular Toeplitz $T$ with $\beta=\min\{m,n\}$, and 
$S_{\alpha}(\eta):=
\begin{bsmallmatrix}
0 & U_{\alpha}(\eta) \\
J_{\alpha}(-i\eta)\overline{U}_{\alpha}(\eta)& 0
\end{bsmallmatrix}$ with $U_{\alpha}(\eta)$
as any solution of $U_{\alpha}(\eta)J_{\alpha}(-\eta^{2})=(J_{\alpha}(i\eta))^{2}U_{\alpha}(\eta)$ for $\alpha\in \{m,n\}$.
\vspace{-1mm}
\item \label{BHHL} If 
$M=L_m(\xi)$ and 
$N=L_n(\xi)$ with 
$\Ima (\xi)>0$ and $\xi^{2}$ nonreal, then 
$Y$ satisfi\-es (\ref{eqAoX2}) if and only if
%
$
Y=R_{m}^{-1}XR_n
$,
in which $X=X_1\oplus \overline{X}_1$ and $X_1$ is an $m$-by-$n$ matrix of the form (\ref{QTY}) for $T$ an $\beta$-by-$\beta$ complex upper triangular Toe\-pli\-tz matrix with $\beta=\min\{m,n\}$, and 
$R_{\alpha}:=P_{\alpha}\oplus P_{\alpha}
$ for $\alpha\in \{m,n\}$.  
\end{enumerate} 
\end{lemma}

The proof of the lemma relies very much on the ideas in \cite{BHH}.

\begin{proof}[Proof of Lemma \ref{lemaBHH}]
The following is a part of Hong's construction of the ca\-no\-nical form under consimilarity \cite[p. 9-10]{Hong90}:
\begin{align*}
H_{\alpha}(\lambda)=P_{\alpha}^{-1}J_{m}(\lambda)\overline{P}_{\alpha},
\,\,\,\,
K_{\alpha}(\mu)=Q_{\alpha}^{-1}\begin{bsmallmatrix}
0  & J_{\alpha}(\mu) \\
-J_{\alpha}(\mu) &  0    
\end{bsmallmatrix}\overline{Q}_{\alpha},
\,\,\,\,
L_{\alpha}(\xi)=R_{\alpha}^{-1}\begin{bsmallmatrix}
0  & J_{\alpha}(\xi) \\
J_{\alpha}(\overline{\xi}) &  0    
\end{bsmallmatrix}\overline{R}_{\alpha},
\end{align*}
in which $\lambda\geq 0$, $\mu>0$, $\xi^{2}\in \mathbb{C}\setminus \mathbb{R}$, and $P_{\alpha}$, $Q_{\alpha}$, $R_{\alpha}$ are as defined in the lemma.

The equation $H_m(\lambda)\overline{Y}=YH_n(\kappa)$ for $\lambda,\kappa\geq 0$ transforms to $J_m(\lambda)\overline{X}=XJ_n(\lambda)$ with $X=P_mYP_n^{-1}$. By setting $X=U+iV$ with real $m$-by-$n$ matrices $U$, $V$, we get $J_m(\lambda)U=UJ_n(\kappa)$ and $-J_m(\lambda)V=VJ_n(\kappa)$. The first equation for $\lambda\neq\kappa$ implies $U=0$, while for $\lambda=\kappa$ we get $U$ upper triangular Toeplitz (see Theorem \ref{lemas}). We write the second equation as $J_m(-\lambda)FV=FVJ_n(\kappa)$ with $F=-1\oplus 1\oplus -1\oplus \cdots$. 
If either $\lambda\neq\kappa$ or $\lambda=\kappa>0$, then $V=0$. When $\lambda=\kappa=0$, then $FV$ is real upper triangular Toeplitz, hence $X$ is complex-alternating upper triangular Toeplitz. This proves (\ref{BHH1}) for $M=H_m(\lambda)$,
$N=H_n(\mu)$ with $\lambda\neq \nu$ and (\ref{BHH2}).

If $V_{\alpha}$ and $S_{\alpha}(\mu)$ are defined as (\ref{BHHK}), it is not difficult to check that
\[
\begin{bsmallmatrix}
0 & J_{\alpha}(\eta) \\
-J_{\alpha}(\eta) & 0
\end{bsmallmatrix}
=
V_{\alpha}^{-1}
\begin{bsmallmatrix}
0 & J_{\alpha}(i\eta) \\
J_{\alpha}(-i\eta) & 0
\end{bsmallmatrix}
\overline{V}_{\alpha}
, \qquad 
S_{\alpha}^{-1}(\eta)\begin{bsmallmatrix}
0 & J_{\alpha}(i\eta) \\
J_{\alpha}(-i\eta) & 0
\end{bsmallmatrix}\overline{S}_{\alpha}(\eta)=\begin{bsmallmatrix}
0 & I_{\alpha}\\
J_{\alpha}(-\eta^{2}) & 0
\end{bsmallmatrix}.
\]
Thus $K_m(\mu)\overline{Y}=YK_n(\nu)$ for $\mu,\nu>0$ transforms to 
\[
J_{m}'(\mu)\overline{X}=XJ_{n}'(\nu),\qquad X=S_m^{-1}(\mu)V_mQ_m Y Q_n^{-1}V_n^{-1}S_n(\nu), \quad 
J_{\alpha}'(\mu):=\begin{bsmallmatrix}
0 & I_{\alpha} \\
J_{\alpha}(-\eta^{2}) & 0
\end{bsmallmatrix}.
\]
Set
$X=\begin{bsmallmatrix}
X_1 & X_2 \\
X_3 & X_4
\end{bsmallmatrix}$:  $\overline{X}_3=X_2J_n(-\nu^{2})$, $J_m(-\mu^{2})\overline{X}_1=X_4J_n(-\nu^{2})$, $J_m(-\mu^{2})\overline{X}_2=X_3$, $\overline{X}_4=X_1$. If 
$\mu=\nu$ we get (\ref{BHHK}), while $\mu\neq \nu$ gives 
(\ref{BHH1}) for $M=K_m(\mu)$,
$N=K_n(\nu)$.

We transform $L_m(\xi)\overline{Y}=YL_n(\zeta)$ for $\Ima (\xi),\Ima (\zeta)>0$ to
\vspace{-1mm}
\[
\begin{bsmallmatrix}
0 & J_{m}(\xi) \\
J_{m}(\overline{\xi}) & 0
\end{bsmallmatrix}\overline{X}
=
X\begin{bsmallmatrix}
0 & J_{n}(\zeta) \\
J_{n\gamma_l}(\overline{\zeta}) & 0
\end{bsmallmatrix}, \qquad   R_m Y R_n^{-1}=X
:=\begin{bsmallmatrix}
X_1 & X_2\\
X_3 & X_4
\end{bsmallmatrix},
\]
where $X_1,X_2,X_3,X_4$ are $m$-by-$n$ matrices.
We have 
\vspace{-1mm}
\begin{align}\label{eqRJ}
&X_{2}J_{n}(\overline{\zeta})=J_{m}(\xi)\overline{X}_{3},\quad 
X_{3}J_{n}(\zeta)=J_{m}(\overline{\xi})\overline{X}_{2},\\
&X_{1}J_{n}(\zeta)=J_{m}(\xi)\overline{X}_{4},\quad 
X_{4}J_{n}(\overline{\zeta})=J_{m}(\overline{\xi})\overline{X}_{1}.\nonumber
\end{align}
By combining the first and the second pair of equations we deduce, respectively, $\overline{X}_{3}(J_{\gamma_l}(\overline{\zeta}))^{2}=(J_{m}(\xi))^{2}\overline{X}_{3}$, $\overline{X}_{2}(J_{n}(\zeta))^{2}=(J_{m}(\overline{\xi}))^{2}\overline{X}_{2}$ and $\overline{X}_{4}(J_{\gamma_l}(\zeta))^{2}=(J_{m}(\xi))^{2}\overline{X}_{4}$, $\overline{X}_{1}(J_{n}(\overline{\zeta}))^{2}=(J_{m}(\overline{\xi}))^{2}\overline{X}_{1}$. Since $\Ima (\xi),\Ima (\zeta)>0$, the first two equations imply $X_3=X_2=0$, while the last two for $\xi\neq \zeta $ yield $X_1=X_4=0$ (thus (\ref{BHH1}) for $M=L_m(\xi)$,
$N=L_n(\zeta)$). Subtracting the third and the last conjugated equation of (\ref{eqRJ}) for $\xi=\zeta$ gives 
$(X_{1}-\overline{X}_4)J_{n}(\xi)=-J_{m}(\xi)(X_1-\overline{X}_{4})$. Hence
$F(X_{1}-\overline{X}_4)J_{n}(\xi)=J_{m}(-\xi)F(X_1-\overline{X}_{4})$, $F=-1\oplus 1\oplus -1\oplus \cdots$, thus we obtain $X_4=\overline{X}_1$. Using (\ref{eqRJ}) then yields that $X_1$
is complex upper triangular Toeplitz, which shows (\ref{BHHL}).
%

Similarly, $K_m(\mu)\overline{Y}=Y L_n(\xi)$ for $\mu>0$, $\xi^{2}\in \mathbb{C}\setminus \mathbb{R} $ reduces to
$\begin{bsmallmatrix}
0 & J_m(\mu) \\
-J_m(\mu) & 0
\end{bsmallmatrix}\overline{Y}=Y\begin{bsmallmatrix}
0 & J_n(\xi) \\
J_n(\overline{\xi}) & 0
\end{bsmallmatrix}$ with $Q_mXR^{-1}_n=Y:=\begin{bsmallmatrix}
X_1 & X_2\\
X_3 & X_4
\end{bsmallmatrix}$ and $X_1,X_2,X_3,X_4$ of size $m\times n$. Thus
\vspace{-1mm}
\begin{align*}
X_{2}J_{n}(\overline{\xi})=J_{m}(\mu)\overline{X}_{3},\quad 
X_{3}J_{n}(\xi)=-J_{m}(\mu)\overline{X}_{2},\\
X_{1}J_{n}(\xi)=J_{m}(\mu)\overline{X}_{4},\quad 
X_{4}J_{n}(\overline{\xi})=-J_{m}(\mu)\overline{X}_{1}.\nonumber
\end{align*}
By combining these equations we get $\overline{X}_{3}(J_n(\overline{\xi}))^{2}=-(J_{m}(\mu))^{2}\overline{X}_{3}$ and $\overline{X}_{4}(J_{n}(\xi))^{2}=-(J_{m}(\mu))^{2}\overline{X}_{4}$, which implies $X_1=X_2=X_3=X_4=0$, hence $X=0$.

Next, 
$H_m(\lambda)\overline{Y}=Y
K_n(\mu)$ for $\lambda\geq 0$, $\mu>0$ reduces to
$J_m(\lambda)\overline{X}=X\begin{bsmallmatrix}
0 & J_n(\mu) \\
-J_n(\mu) & 0
\end{bsmallmatrix}$, whe\-re $P_mYQ^{-1}_n=X:=\begin{bsmallmatrix}
X_1 & X_2
\end{bsmallmatrix}$ with $m$-by-$n$ matrices $X_1,X_2$.
We get $J_m(\lambda)\overline{X}_1=-X_2J_n(\mu)$, $J_m(\lambda)\overline{X}_2=X_1J_n(\mu)$, thus $(J_m(\lambda))^{2}X_1=-J_m(\lambda)\overline{X}_2J_n(\mu)=-X_1(J_n(\mu))^{2}$. 
It yields 
$S^{-1}J_m(\lambda^{2})SX_1=-X_1T^{-1}FJ_n(-\mu^{2})F^{-1}T$
and for some nonsingular $S$, $T$ and $F=-1\oplus 1\oplus -1\oplus \cdots$. 
Since $\lambda^{2}\geq 0>-\mu^{2}$, Theorem \ref{lemas} implies $SX_1T^{-1}F=0$ with $X_1=0$ (hence $X_2=0$), and therefore $X=0$.

Further, 
$H_m(\lambda)\overline{Y}=Y
L_n(\xi)$ 
yields 
$J_m(\lambda)\overline{X}=X\begin{bsmallmatrix}
0 & J_n(\xi) \\
J_n(\overline{\xi}) & 0
\end{bsmallmatrix}$ with $P_mYR^{-1}_n=X:=\begin{bsmallmatrix}
X_1 & X_2
\end{bsmallmatrix}$ for some $m$-by-$ n$ matrices $X_1,X_2$.
We obtain equations $J_m(\lambda)\overline{X}_1=X_2J_n(\overline{\xi})$ and $J_m(\lambda)\overline{X}_2=X_1J_n(\xi)$, therefore $(J_m(\lambda))^{2}X_1=J_m(\lambda)\overline{X}_2J_n(\xi)=X_1(J_n(\xi))^{2}$. If $\lambda\geq 0$ and $\xi^{2}$ is nonreal, Theorem \ref{lemas} yields $X_1=X_2=0$), thus $X=0$.

Since $H_m(\lambda)$, $K_n(\mu)$, $L_n(\xi)$ are Hermitian, by conjugating and transposing equations $K_n(\mu)\overline{Y}=YH_m(\lambda)$, $L_n(\xi)\overline{Y}=YH_m(\lambda)$, $L_n(\xi)\overline{Y}=YK_m(\mu)$ we obtain $Y^{T}K_n(\mu)=H_m(\lambda)\overline{Y}^{T}$, $L_n(\xi)Y^{T}=H_m(\lambda)\overline{Y}^{T}$, $Y^{T}L_m(\xi)=K_n(\mu)\overline{Y}^{T}$, respectively. These equations have already been solved with solution $Y=0$. This concludes 
(\ref{BHH1}).
\end{proof}

\begin{remark}
The form of a solution of (\ref{eqAoX2}) for $M=L_m(\xi)$, $N=L_n(\xi)$ with $\xi^{2}\in \mathbb{C}\setminus \mathbb{R}$ in \cite{BHH} is not suited for our application in the proof of Theorem \ref{stabz}; the usage of $\begin{bsmallmatrix}
0 & J_m(\xi)\\
J_m(\overline{\xi}) & 0
\end{bsmallmatrix}$ instead of 
$\begin{bsmallmatrix}
0 & I_m \\
J_m(\xi^{2}) & 0
\end{bsmallmatrix}$ in the proof of Lemma \ref{lemaBHH} is essential.
\end{remark}

We proceed with a technical lemma based on the idea from the paper by Lin, Mehrmann and Xu \cite[Sec. 3.1]{Lin} (see also \cite[Sec. 2]{TSOS}). It enables us to transform a block matrix with (complex-alternating) upper triangular Toeplitz blocks to a block (complex-alternating) upper triangular Toeplitz matrix. 
Set
\begin{equation}\label{perS}
\Omega_{\alpha,m}:=\left[e_1\;e_{\alpha+1}\;\ldots\;e_{(m-1)\alpha+1}\;e_2\;e_{\alpha+2}\;\ldots\;e_{(m-1)\alpha+2}\;\ldots\;e_{\alpha}\;e_{2\alpha}\;\ldots\;e_{\alpha m}\right],
\end{equation}
where $e_1,e_2,\ldots,e_{\alpha m}$ is the standard orthonormal basis in $\mathbb{C}^{\alpha m}$.
Multiplication with $\Omega_{\alpha,m}$ from the right (with $\Omega_{\alpha,m}^{T}$ from the left) puts the $k$-th, the $(\alpha+k)$-th, \ldots, the $((m-1)\alpha+k)$-th column (row) together for all $ k \in \{1,\ldots, \alpha\}$. 
For example,
\vspace{-1mm}
\[
\Omega_{3,2}^{T}\begin{bmatrix}[cc|cc|cc]
a_1 & b_1 & a_2 & b_2 & a_3 & b_3 \\
0   & a_1 & 0   & a_2 & 0   & a_3\\
0   & 0   & 0   & 0   & 0   &  0 \\
\hline
a_4 & b_4 & a_5 & b_5 & a_6 & b_6 \\
0   & a_4 & 0   & a_5 & 0   & a_6\\
0   & 0   & 0   & 0   & 0   &  0 
\end{bmatrix}\Omega_{2,3}
=
\begin{bmatrix}[ccc|ccc]
a_1 & a_2 & a_3 & b_1 & b_2 & b_3 \\
a_4 & a_5 & a_6   & b_4 & b_5   & b_6\\
\hline
0   & 0   & 0   & a_1   & a_2   &  a_3 \\
0 &  0 &   0 & a_4 & a_5 & a_6 \\
\hline
0   & 0 & 0   & 0 & 0   & 0\\
0   & 0   & 0   & 0   & 0   &  0 
\end{bmatrix}.
\]
%
%
Similarly, multiplication with 
the following matrix from the right 
puts the $k$-th, the $(2\alpha+k)$-th,\ldots,the $((2m-1)\alpha+k)$-th column (row) together:
\vspace{-1mm}
\small
\begin{align}\label{perK}
\Omega_{\alpha,m}':=\Big[
& e_1\;e_{2\alpha+1}\;\ldots\;e_{2(m-1)\alpha+1}\;
e_{\alpha+1}\;e_{3\alpha+1}\;\ldots\;e_{(2m-1)\alpha+1}\;
e_2\;e_{2\alpha+2}\;\ldots\;e_{2(m-1)\alpha+2}\;
\nonumber
\\
& e_{\alpha+2}\;e_{3\alpha+2}\;\ldots\;e_{(2m-1)\alpha+2}\;\ldots\;
\ldots\; e_{\alpha}\;e_{3\alpha}\;\ldots\;e_{\alpha (2m-1)}\;
e_{2\alpha}\;e_{4\alpha}\;\ldots\;e_{\alpha (2m)} 
\Big].
\end{align}
\normalsize
It is then immediate:

\begin{lemma}\label{lemaP}
Suppose 
$X=[X_{rs}]_{r,s=1}^{N}$ such that each block $X_{rs}=[(X_{rs})_{jk}]_{j,k=1}^{m_r,m_s}$ is an $m_r$-by-$m_s$ block matrix with blocks of the same size, and let $\alpha_{1}>\ldots >\alpha_{N}$ with $b_{rs}:=\{\alpha_r,\alpha_s\}$. Also, set $\Omega:=\bigoplus_{r=1}^{N}\Omega_{\alpha_r,m_r}$ and $\Omega':=\bigoplus_{r=1}^{N}\Omega_{\alpha_r,m_r}'$.
\begin{enumerate}
\item \label{lemaP1}
Assume that each $X_{rs}$ 
con\-sists of blocks of size $\alpha_r\times \alpha_s$ 
and such that 
\vspace{-1mm}
\[
(X_{rs})_{jk}=
\left\{
\begin{array}{ll}
[0\quad T_{jk}^{rs}], & \alpha_{r}<\alpha_{s}\\
\begin{bmatrix}
T_{jk}^{rs}\\
0
\end{bmatrix}, & \alpha_{r}>\alpha_{s}\\
T_{jk}^{rs},& \alpha_{r}=\alpha_{s}
\end{array}\right., \qquad
\begin{array}{l}
T_{jk}^{rs}=T(a_{0,jk}^{rs},a_{1,jk}^{rs},\ldots,a_{b_{rs}-1,jk}^{rs})\\
\\
(\textrm{or }T_{jk}^{rs}=T_c(a_{0,jk}^{rs},a_{1,jk}^{rs},\ldots,a_{b_{rs}-1,jk}^{rs}))
\end{array}
\]
for $j\in \{1,\ldots, m_r\}$, $k\in \{1,\ldots, m_s\}$, $a_{n,jk}^{rs}\in \mathbb{C}$,
and set $ A_{n}^{rs}:=[a_{n,jk}^{rs}]_{j,k=1}^{m_r,m_s}$.
Then 
\vspace{-1mm}
\begin{align}\label{0T02}
&\mathcal{X}:=\Omega^{T}X\Omega, \quad \mathcal{X}=[\mathcal{X}_{rs}]_{r,s=1}^{N},\qquad
\mathcal{X}_{rs}=
\left\{
\begin{array}{ll}
[0\quad \mathcal{T}_{rs}], & 
\alpha_r<\alpha_s\\
\begin{bmatrix}
\mathcal{T}_{rs}\\
0
\end{bmatrix}, & \alpha_r>\alpha_s\\
\mathcal{T}_{rs},& \alpha_r=\alpha_s
\end{array}\right.,
\end{align}
%
with $\mathcal{X}_{rs}$ of size $\alpha_r \times \alpha_s$ and
$\mathcal{T}_{rs}\hspace{-0.3mm}=\hspace{-0.3mm}
T(A_0^{rs},\ldots,A_{b_{rs}-1}^{rs})$ ($\mathcal{T}_{rs}\hspace{-0.3mm}=\hspace{-0.3mm}
T_c(A_0^{rs},\ldots,A_{b_{rs}-1}^{rs})$).
%
\item \label{lemaP2}
Let each $X_{rs}$ consists of four blocks of size $\alpha_r\times \alpha_s$, and such that:
\vspace{-1mm}
\begin{align*}
&(X_{rs})_{jk}=\begin{bmatrix}
\tau_{jk}^{rs} & \sigma_{jk}^{rs}\\
J_{\alpha_r}(\eta)\overline{\sigma}_{jk}^{rs} & \overline{\tau}_{jk}^{rs}
\end{bmatrix}, \qquad  j\in \{1,\ldots m_r\}, \quad k\in \{1,\ldots m_s\},\quad \rho\in \mathbb{C},\\
&\tau_{jk}^{rs}=
\left\{
\begin{array}{ll}
\begin{bsmallmatrix}
0 & T_{jk}^{rs}
\end{bsmallmatrix}, & \alpha_{r}<\alpha_{s}\\
\begin{bsmallmatrix}
T_{jk}^{rs}\\
0
\end{bsmallmatrix}, & \alpha_{r}>\alpha_{s}\\
T_{jk}^{rs},& \alpha_{r}=\alpha_{s}
\end{array}\right., 
\quad
\sigma_{jk}^{rs}=
\left\{
\begin{array}{ll}
\begin{bsmallmatrix}
0 & S_{jk}^{rs}
\end{bsmallmatrix}, & \alpha_{r}<\alpha_{s}\\
\begin{bsmallmatrix}
S_{jk}^{rs}\\
0
\end{bsmallmatrix}, & \alpha_{r}>\alpha_{s}\\
S_{jk}^{rs},& \alpha_{r}=\alpha_{s}
\end{array}\right.,\\
& T_{jk}^{rs}=T(v_{0,jk}^{rs},\ldots,v_{b_{rs}-1,jk}^{rs}), \,\,\, S_{jk}^{rs}=T(w_{0,jk}^{rs},\ldots,w_{b_{rs}-1,jk}^{rs}), \quad\textrm{all }v_{n,jk}^{rs},w_{n,jk*}^{rs}\in \mathbb{C}.
\end{align*}
%
Set $W_{-1}^{rs}:=0$ and $V_n^{rs}:=[v_{n,jk}^{rs}]_{j,k=1}^{m_r,m_s}$, 
$W_n^{rs}:=[w_{n,jk}^{rs}]_{j,k=1}^{m_r,m_s}$ with
$A_n^{rs}:=\begin{bsmallmatrix}
V_n^{rs} & W_n^{rs}\\
\rho \overline{W}_n^{rs}+\overline{W}_{n-1}^{rs} & \overline{V}_n^{rs}
\end{bsmallmatrix}$ for $n\in \{0,\ldots,b_{rs}-1\} $.
Then 
\begin{align*}
\mathcal{X}':=(\Omega')^{T}X\Omega'
\end{align*}
is of the form (\ref{0T02}) with $\mathcal{T}_{rs}=T(A_0^{rs},\ldots,A_{b_{rs}-1}^{rs})$.

\quad
Furthermore, if all $W_n^{rs}=0$,
then there exists a permutation matrix $\Omega_0$ such that $\Omega_0^{T}X\Omega_0 =\mathcal{V}\oplus\overline{\mathcal{V}}$ with $\mathcal{V}$ of the form (\ref{0T02}) for $\mathcal{T}_{rs}=T(V_0^{rs},\ldots,V_{b_{rs}-1}^{rs})$.
\end{enumerate}
\end{lemma}

The following proposition describes the (nonsingular) solutions of (\ref{eqAoXXB}). 

\begin{proposition} \label{resAoXXA}
\begin{enumerate}[label={\arabic*.},ref={\arabic*},
leftmargin=20pt]
\item \label{resAoXXA1}
Let $\rho_1,\ldots,\rho_n\in \mathbb{C}$ be all distinct and let 
$
\mathcal{H}=\bigoplus_{j=1}^{n}\mathcal{H}_j$,
in which $\mathcal{H}_j$ is a direct sum whose summands 
are either of the form (\ref{Hmz}) or (\ref{KLmz}), and such that they  
correspond to the eigenvalue $\rho_j$ of $\mathcal{H}\overline{\mathcal{H}}$.  
Then the solution of $\mathcal{H}\overline{Y}=Y\mathcal{H}$ is of the form $Y=\bigoplus_{j=1}^{n} Y_j$ with $Y_j$ as a solution of $\mathcal{H}_j\overline{Y}_j=Y_j\mathcal{H}_j$.

\item \label{resAoXXA2}
For $\mu=(m_1,\ldots,m_N)$, $\alpha=(\alpha_1,\ldots,\alpha_N)$ let $\mathbb{T}^{\alpha,\mu}$, $\mathbb{T}^{\alpha,\mu}$ and $\mathbb{T}_c^{\alpha,\mu}$ consist of matrices as described in (\ref{0T0}),
and let $\mathcal{H}=\mathcal{H}^{\varepsilon}$ be as in Theorem \ref{stabw} for all $\varepsilon_{r,j}=1$ ($\mathcal{H}\overline{\mathcal{H}}$ has precisely one eigenvalue $\rho$).
The nonsingular solutions of $\mathcal{H}\overline{Y}=Y\mathcal{H}$ form a group conjugate to $\mathbb{T}^{\alpha,\mu}\oplus\overline{\mathbb{T}}^{\alpha,\mu}$ for $\rho\in \mathbb{C}\setminus \mathbb{R}$, conjugate to $\mathbb{T}_c^{\alpha,\mu}$ for $\rho=0$, conjugate to the subgroup of all real matrices in $\mathbb{T}^{\alpha,\mu}$ for $\rho>0$, and conjugate to the subgroup of all matrices in $\mathbb{T}^{\alpha,2\mu}$ of the form (\ref{0T0}) with (\ref{Trsvw0}) for $\rho<0$.
\end{enumerate}
\end{proposition}

\begin{proof}
Suppose 
$\mathcal{H}=\bigoplus_{j}^{} M_{j}$
with all 
$M_j$ either of the form (\ref{Hmz}) or of the form (\ref{KLmz}).
The equation $\mathcal{H}\overline{Y}=Y\mathcal{H}$ is then equivalent to a system of equations:
\begin{equation}\label{eqTsys}
M_{j}\overline{Y}_{jk}=Y_{jk}M_{k}, \quad j,k=1,2,\ldots, \qquad 
Y:=[Y_{jk}]_{j,k}^{},
\end{equation}
in which $Y$ is partitioned conformally to $\mathcal{H}$. Lemma \ref{lemaBHH} (\ref{BHH1}) implies (\ref{resAoXXA1}.

Next, let all $M_j\overline{M}_j$ have the same eigenvalue $\rho$.
In view of Lemma \ref{lemaBHH} there exist nonsingular matrices $U_j$ so that any solution $Y$ of (\ref{eqTsys}) is of the form 
\vspace{-1mm}
\begin{equation*}
Y=U^{-1}XU \quad  (Y_{jk}=U_{j}^{-1}X_{jk}U_{k}^{-1}); \qquad
X:=[X_{jk}]_{j,k}^{}, U:=\oplus \textrm{}_{j}^{} U_{j},
\end{equation*}
%
where all $X_{jk}$ are of the form (\ref{QTY}) with real (complex-alternating) upper triangular Toeplitz $T$ for $\rho>0$ (for $\rho=0$), or of the form  (\ref{QTC}) with upper triangular Toe\-p\-litz $X_1$, $X_2$ (and $X_2=0$) for $\rho<0$ (for $\rho\in \mathbb{C}\setminus \mathbb{R}$).
Lemma \ref{lemaP} gives (\ref{resAoXXA2}).
\end{proof}

We observe 
the group structures of 
$\mathbb{T}^{\alpha,\mu}$, $\mathbb{T}_c^{\alpha,\mu}$.
The claim for $\mathbb{T}^{\alpha,\mu}$ 
coincid\-es with \cite[Lemma 2.2]{TSOS} and its proof is based on ideas from 
\cite[Example 6.49]{Milne} describing upper unitriangular matrices; 
it works mutatis mutandis for $\mathbb{T}_c^{\alpha,\mu}$.

\begin{lemma}\label{lemanilpo}
Let $\mathbb{T}^{\alpha,\mu}$ and $\mathbb{T}_c^{\alpha,\mu}$ consist of matrices defined in (\ref{0T0}). 
Then 
$\mathbb{T}^{\alpha,\mu}=\mathbb{D}\ltimes \mathbb{U}$ and $\mathbb{T}_c^{\alpha,\mu}=\mathbb{D}_c\ltimes \mathbb{U}_c$ are semidirect products of subgroups, where $\mathbb{D}\subset \mathbb{T}^{\alpha,\mu}$, $\mathbb{D}_c\subset \mathbb{T}_c^{\alpha,\mu}$ contain nonsingular block diagonal matrices, and $\mathbb{U}\subset \mathbb{T}^{\alpha,\mu}$, $\mathbb{U}_c\subset \mathbb{T}_c^{\alpha,\mu}$ are normal subgroups consisting of upper (complex-alternating) unitriangular Toeplitz diagonal blocks. 
Moreover, $\mathbb{U}$ and $\mathbb{U}_c$ are unipotent of order $\leq \alpha_1-1$.
\end{lemma}

\section{Certain block matrix equation}\label{cereq}

Let $\alpha_{1}>\alpha_{2}>\ldots >\alpha_{N}$ and suppose that
we are given nonsingular matrices
\vspace{-1mm}
\small
\begin{align}\label{BBF}
&\mathcal{B}=\bigoplus_{r=1}^{N}T\big(B_0^{r},B_1^{r},\ldots,B_{\alpha_r-1}^{r}\big),\quad \mathcal{C}=\bigoplus_{r=1}^{N}T\big(C_0^{r},C_1^{r},\ldots,C_{\alpha_r-1}^{r}\big),\quad \mathcal{F}=\bigoplus_{r=1}^{N}E_{\alpha_r}(I_{m_r}),
%
\vspace{-1mm}
\end{align}
\normalsize
%
with symmetric 
$B_n^{r},C_n^{r}\in  \mathbb{C}^{m_r \times m_r}$
and
$E_{\beta}(I_{m}):=
\begin{bsmallmatrix}
 0                 &      & I_{m}\\
            &   \iddots     &  \\
I_{m}            &           &  0\\
\end{bsmallmatrix}
$
is an $\beta\textrm{-by-} \beta$ block matrix with $I_m$ on the anti-diagonal and zero matrices otherwise.
%
We find all $\mathcal{X}$ in $ \mathbb{T}^{\alpha,\mu}$ or $ \mathbb{T}_c^{\alpha,\mu}$ for $\alpha=(\alpha_1,\ldots,\alpha_N)$, $\mu=(m_1,\ldots,m_N)$ (see (\ref{0T0}))
that solve 
%
\vspace{-1mm}
\begin{equation}\label{eqFYFIY}
\mathcal{C}=\mathcal{F}\mathcal{X}^{T}\mathcal{F}\mathcal{B} \mathcal{X};
\end{equation}
%
this is essential to prove Theorem \ref{stabw} and Theorem \ref{stabz}. The observation
%
%
\vspace{-1mm}
\[
(\mathcal{F}\mathcal{X}^{T}\mathcal{F}\mathcal{B} \mathcal{X})^{T}=\mathcal{X}^{T}\mathcal{B}^{T}\mathcal{F}\mathcal{X}\mathcal{F} =\mathcal{F}\mathcal{F}\mathcal{X}^{T}\mathcal{F}(\mathcal{F}\mathcal{B}^{T}\mathcal{F})\mathcal{X}\mathcal{F}=\mathcal{F}(\mathcal{F}\mathcal{X}^{T}\mathcal{F}\mathcal{B}\mathcal{X})\mathcal{F}
\]
shows that for $r\neq s$ we have $(\mathcal{F}\mathcal{X}^{T}\mathcal{F}\mathcal{B}\mathcal{X})_{rs}=0$ if and only if $(\mathcal{F}\mathcal{X}^{T}\mathcal{F}\mathcal{B}\mathcal{X})_{sr}=0$.
When comparing the left-hand side with the right-hand side of (\ref{eqFYFIY}) blockwise, it thus suffices to observe the upper triangular parts of $\mathcal{F}X^{T}\mathcal{F}\mathcal{B} X$ and $\mathcal{C}$.  
Since $(\mathcal{F}\mathcal{X}^{T}\mathcal{F}\mathcal{B}\mathcal{X})_{rs}$ and $\mathcal{C}_{rs}$ are rectangular upper triangular Toeplitz of the same form,
it is enough to compare their first rows.
By simplifying the notation with $\mathcal{Y}:=\mathcal{B}\mathcal{X}$ and $\widetilde{\mathcal{X}}:=\mathcal{F}X^{T}\mathcal{F}$, we obtain 
the entry in the $j$-th column and in the first row of $(\mathcal{F}\mathcal{X}^{T}\mathcal{F}\mathcal{B}\mathcal{X})_{rs}=(\widetilde{\mathcal{X}}\mathcal{Y})_{rs}$ 
by multiplying the first rows of blocks $\widetilde{\mathcal{X}}_{r1},\ldots, \widetilde{\mathcal{X}}_{rN}$ with the $j$-th columns of blocks $\mathcal{Y}_{1s},\ldots, \mathcal{Y}_{Ns}$, respectively, and then adding them.
Hence (\ref{eqFYFIY}) redudes to:
\vspace{-1mm}
\begin{align}
\label{f1}
(\mathcal{C}_{r(r+p)})_{1j}= & (\widetilde{\mathcal{X}}_{rr})_{(1)}(\mathcal{Y}_{r(r+p)})^{(j)}+\sum_{k=r+1}^{N}(\widetilde{\mathcal{X}}_{rk})_{(1)}(\mathcal{Y}_{k(r+p)})^{(j)}\\
\vspace{-1mm}
  &+\sum_{k=1}^{r-1}(\widetilde{\mathcal{X}}_{rk})_{(1)}(\mathcal{Y}_{k(r+p)})^{(j)},
 \quad 1\leq j\leq \alpha_{r+p},\quad 0\leq p \leq N-r.
  \nonumber
\end{align}
It turns out to be important to consider equations (\ref{f1}) in an appropriate order.
The following lemmas provide this computation 
in detail.

\begin{lemma}\label{EqT}
\hspace{-1mm}Let $\mathcal{B},\mathcal{C}$ as in (\ref{BBF}) be given.\hspace{-1mm} Then the 
dimension of the space of soluti\-ons
of (\ref{eqFYFIY}) that are of the form $\mathcal{X}\hspace{-1mm}=\hspace{-0.5mm}[\mathcal{X}_{rs}]_{r,s=1}^{N}$ \hspace{-0.5mm}(partitioned conformally to $\mathcal{B},\mathcal{C}$)
with 
%
\vspace{-1mm}
\begin{equation}\label{EqTX}
\mathcal{X}_{rs}=
\left\{
\begin{array}{ll}
\hspace{-1mm}[0\quad \mathcal{T}_{rs}], & \alpha_r<\alpha_s\\
\hspace{-1mm}\begin{bmatrix}
\mathcal{T}_{rs}\\
0
\end{bmatrix}, & \alpha_r>\alpha_s\\
\hspace{-1mm}\mathcal{T}_{rs},& \alpha_r=\alpha_s
\end{array}\right., \quad 
\begin{array}{l}
(\alpha_{1}>\alpha_{2}>\ldots >\alpha_{N}),\\
b_{rs}:=\min\{\alpha_s,\alpha_r\},\\
\mathcal{T}_{rs}=T\big(A_0^{rs},\ldots,A_{b_{rs}-1}^{rs}\big),\,\,A_j^{rs}\in \mathbb{C}^{m_r\times m_s}
\end{array}
\end{equation}
is $\sum_{r=1}^{N}m_r(\tfrac{m_r-1}{2}\alpha_r+\sum_{s=1}^{r-1}\alpha_s m_s)$, and each solution
satisfies the following properties:
\begin{enumerate}[label={(\alph*)},ref={\alph*},
leftmargin=13pt,itemindent=3pt]
%
\item \label{EqT2} 
Each $A_0^{rr}$ is a solution of the equation $C_0^{r}=(A_0^{rr})^{T}B_0^{r}A_0^{rr}$.
If $N\geq 2$ matrices $A_j^{rs}$ for $j\in \{0,\ldots,\alpha_{r}-1\}$, $r,s\in \{1,\ldots,N\}$ with $r>s$ can be taken freely. 
\item \label{EqT3} 
Assuming (\ref{EqT2}) and choosing
matrices 
$Z_j^{r}=-Z_j^{r}\in \mathbb{C}^{m_r\times m_r}$ for $r\in \{1,\ldots,N\}$, $j\in \{1,\ldots,\alpha_r-1\}$ freely, 
the remaining entries of $\mathcal{X}$ are computed 
as follows:

\vspace{-3mm}
\hspace{-13mm}

\begin{algorithmic}
\State $\Psi_n^{krs}:=\sum_{i=0}^{n} \sum_{l=0}^{n-i}(A_i^{kr})^{T} B_{n-i-l}^{k}A_l^{ks}$
\State $\widetilde{\Psi}_n^{krs}:=\sum_{i=1}^{n} \sum_{l=0}^{n-i}(A_i^{kr})^{T} B_{n-i-l}^{k}A_l^{ks}+ \sum_{l=0}^{n-1}(A_0^{kr})^{T} B_{n-l}^{k}A_l^{ks}$
\hspace{-3mm}\For {$j=0:\alpha_1-1$}
    \If {$r\in \{1,\ldots,N\}$, $ j\in \{1,\ldots,\alpha_r-1\}$}
    \State \hspace{-5mm} $A_j^{rr}=\frac{1}{2}A_0^{rr}-\frac{1}{2}A_0^{rr}(C_0^{r})^{-1}(Z_j^{r}+\widetilde{\Psi}_j^{rrr}
    +\sum_{k=1}^{r-1}\Psi_{j-\alpha_k+\alpha_r}^{krr}+\sum_{k=r+1}^{N}\Psi_{j-\alpha_{r}+\alpha_k}^{krr})$
    \EndIf
    \For {$p=1:N-1$}
        \If {$r\in \{1,\ldots,N\}$, $j\leq \alpha_{r+p}-1$, $r+p\leq N$}
        \small
        \State \hspace{-5mm} $A_j^{r(r+p)}=-A_0^{r(r+p)}(C_0^{r})^{-1}\big(
        (A_{j}^{rr})^{T} B_{0}^{r}A_0^{r(r+p)}+\widetilde{\Psi}_j^{rr(r+p)}+\sum_{k=1}^{r-1}\Psi_{j-\alpha_k+\alpha_r}^{kr(r+p)}$\\
        \qquad \qquad \qquad \qquad \qquad \qquad \qquad \quad  $+ \sum_{k=r+1}^{r+p}\Psi_{j}^{kr(r+p)}+\sum_{k=r+p+1}^{N}\Psi_{j-\alpha_{r+p}+\alpha_k}^{kr(r+p)}\big)$   
\normalsize
        \EndIf
    \EndFor
\EndFor
\end{algorithmic}
For simplicity, we define $\sum_{j=l}^{n}a_j=0$ if $l>n$, and it is understood that the inner loop (i.e. for p =1 : N-1) is not performed for $N=1$. 

\vspace{-1mm}

\item \label{EqTII} 
\begin{enumerate}[label={(\roman*)},ref={\roman*},
leftmargin=15pt
]
\item \label{EqTIIi}
If $\mathcal{B},\mathcal{C}$ are real, then $\mathcal{X}$ is real 
if and only if 
the following statements hold
\begin{itemize}
\item Matrices $B_0^{r}$ and $C_0^{r}$ in (\ref{BBF}) 
have the same inertia
for all $r\in \{1,\ldots,N\}$.
\item All matrices $A_0^{rr}$, matrices $A_j^{rs}$ with $r>s$, $j\in \{0,\ldots,\alpha_{r}-1\}$,
and $Z_j^{r}$ for $j\in \{1,\ldots,\alpha_{r}-1\}$ in 
(\ref{EqT2}) and (\ref{EqT3}) 
are chosen real.
\end{itemize}

\item \label{EqTIII} 
For any $r\in \{1,\ldots,N\}$, $n\in \{1,\ldots,b_{rs}-1\}$ assume in (\ref{BBF}) that $m_r=2m_r'$ and
%
\vspace{-2mm}
\begin{align}\label{aass}
&B_n^{r}=u_n^{r} K_r+u_{n-1}^{r}L_r, \quad 
K_{r}:=
-\mu^{2}I_{m_r'}
\oplus I_{m_r'}, \quad
L_{r}:=
I_{m_r'} 
\oplus 0,\quad \mu>0,\\
&C_n^{r}=v_n^{r} K_r+v_{n-1}^{r}L_r,\quad u_0,v_0,\ldots,u_{b_{rs}-1},v_{b_{rs}-1}\in \mathbb{R},\,\,u_0,v_0\neq 0,\,\,u_{-1}=v_{-1}=0. \nonumber
\end{align}
Then there are $V_j^{rs}, W_j^{rs}\in \mathbb{C}^{m_r'\times m_r'}$ for $j\in \{0,\ldots,b_{rs}-1\}$ and such that 
\vspace{-2mm}
\begin{align}\label{AVW} 
&A_0^{rs}=
\begin{bsmallmatrix}
V_0^{rs} &  W_0^{rs} \\
-\mu^{2}\overline{W}_0^{rs} &  \overline{V}_0^{rs} 
\end{bsmallmatrix}, \quad
A_n^{rs}=
\begin{bsmallmatrix}
V_n^{rs} &  W_n^{rs} \\
-\mu^{2}\overline{W}_n^{rs}+\overline{W}_{n-1}^{rs} &  \overline{V}_n^{rs} 
\end{bsmallmatrix}, \quad
n\in \{1,\ldots, b_{rs}-1\},
\end{align}
precisely when 
$A_0^{rs},Z_j^{r}$ in (\ref{EqT2}), (\ref{EqT3}) 
are of the form $\begin{bsmallmatrix}
V^{} &  W^{} \\
-\mu^{2}\overline{W}^{} &  \overline{V}^{} 
\end{bsmallmatrix}$, $V^{}, W^{}\in \mathbb{C}^{m_r'\times m_r'}$.
\end{enumerate}

\end{enumerate}
\end{lemma}


Lemma \ref{EqT} (\ref{EqT2}), (\ref{EqT3}), (\ref{EqTII}) (\ref{EqTIIi}) coincides with \cite[Lemma 3.1]{TSOS}; we apologize for minor errors in formulas providing $A_j^{rr}$ and $A_j^{r(r+p)}$ in \cite[Lemma 3.1 (b)]{TSOS}. Thus we 
only prove (\ref{EqTII}) (\ref{EqTIII}), in which solutions are of a special form, which makes the analysis considerably more involved.

\begin{lemma}\label{EqTca}
Let $\mathcal{B}$, $\mathcal{C}$ as in (\ref{BBF}) and real be given.
Then the solution of (\ref{eqFYFIY}) that is of the form
$\mathcal{X}=[\mathcal{X}_{rs}]_{r,s=1}^{N}$ (partitioned conformally to $\mathcal{B},\mathcal{C}$)
%
with
\vspace{-1mm}
\begin{equation}\label{EqTXca}
\mathcal{X}_{rs}=
\left\{
\begin{array}{ll}
[0\quad \mathcal{T}_{rs}], & \alpha_r<\alpha_s\\
\begin{bmatrix}
\mathcal{T}_{rs}\\
0
\end{bmatrix}, & \alpha_r>\alpha_s\\
\mathcal{T}_{rs},& \alpha_r=\alpha_s
\end{array}\right., \quad 
\begin{array}{l}
(\alpha_{1}>\alpha_{2}>\ldots >\alpha_{N}),\\
b_{rs}:=\min\{\alpha_s,\alpha_r\}\\
\mathcal{T}_{rs}=T_c\big(A_0^{rs},\ldots,A_{b_{rs}-1}^{rs}\big),\,\,\, A_j^{rs}\in \mathbb{C}^{m_r\times m_s},
\end{array}
\end{equation}
exists if and only if the following condition holds:
\begin{equation}\label{cond}
B_0^{r} \textrm{ and } C_0^{r} \textrm{ have the same inertia for all } r\in \{1,\ldots,N\} \textrm{ such that } \alpha_r \textrm{ is even}.
\end{equation}
If (\ref{cond}) is fulfilled, then the real dimension of the space of solutions is 
\vspace{-1mm}
\small
\begin{align*}
\sum_{r=1}^{N}\big(\alpha_r m_r^{2} +2\sum_{s=1}^{r-1}\alpha_s m_r m_s\big) - \sum_{\alpha_r \textrm{ even}}\tfrac{\alpha_r}{2}m_r-\sum_{\alpha_r \textrm{ odd}}\tfrac{\alpha_r+1}{2} m_r.
\vspace{-1mm}
\end{align*}
\normalsize
Furthermore, such solutions
satisfy the following properties:
%
\begin{enumerate}[label={(\alph*)},ref={\alph*},
leftmargin=13pt,itemindent=3pt]
\item \label{EqT2ca} 
Each $A_0^{rr}$ with $\alpha_r$ odd is a solution of $C_0^{r}=(A_0^{rr})^{T}B_0^{r}A_0^{rr}$, while $A_0^{rr}$ for $\alpha_r$ even is a solution of $C_0^{r}=(A_0^{rr})^{*}B_0^{r}A_0^{rr}$.
If $N\geq 2$ the entries of $A_j^{rs}$ for $j\in \{0,\ldots,\alpha_{r}-1\}$ and $r,s\in \{1,\ldots,N\}$ with $r>s$ can be taken as free variables.
\vspace{-1mm}
\item \label{EqT3ca} 
Assuming (\ref{EqT2}) and choosing all 
$m_r$-by-$m_r$ matrices 
\small 
$Z_j^{r}=
\left\{
\begin{array}{ll}
\hspace{-1mm}-(Z_j^{r})^{T}, & \hspace{-1mm} j-\alpha_r \textrm{ odd }\\
\hspace{-1mm}-(Z_j^{r})^{*}, & \hspace{-1mm} j-\alpha_r \textrm{ even }
\end{array}
\right.$
\normalsize
for 
$j\in \{1,\ldots,\alpha_r-1\}$ freely, the remaining entries of $\mathcal{X}$ are computed 
as follows:
\end{enumerate}

\vspace{-5mm}
\hspace{-12mm}
\begin{algorithmic}
\State 
\small
$
\mathcal{A}_n^{kr}:=
\left\{
\begin{array}{ll}
\hspace{-0.5mm} 
\begin{bsmallmatrix}
(A_0^{kr})^{T} &(\overline{A}_1^{kr})^{T} &\ldots & (\overline{A}_{n-1}^{kr})^{T} & (A_{n}^{kr})^{T}
\end{bsmallmatrix}, & \hspace{-0.5mm} n \textrm{ even}\\ 
\hspace{-0.5mm}
\begin{bsmallmatrix}
(A_0^{kr})^{T} & (\overline{A}_1^{kr})^{T} & \ldots & (A_{n-1}^{kr})^{T} & (\overline{A}_{n}^{kr})^{T}
\end{bsmallmatrix}, & \hspace{-0.5mm} n \textrm{ odd}
\end{array}
\right.\hspace{-0.5mm}; \quad
\mathcal{R}_{n}^{k}:=
\left\{
\begin{array}{ll}
\hspace{-0.5mm}
\begin{bsmallmatrix}
B_{n}^{k} & B_{n-1}^{k} &\ldots & B_{1}^{k}
\end{bsmallmatrix}, & \hspace{-0.5mm} n\neq 0\\
\hspace{-0.5mm} 
0, & \hspace{-0.5mm} n=0
\end{array}
\right.\hspace{-0.5mm},
$
\State
$\phi_{n}^{ks}:=
\left\{\begin{array}{ll}
\hspace{-0.5mm} \mathcal{R}_n^{k}(\mathcal{A}_{n}^{ks})^{T} & \hspace{-0.5mm} n \textrm{ even}\\
\hspace{-0.5mm} \mathcal{R}_n^{k}(\overline{\mathcal{A}}_n^{ks})^{T} & \hspace{-0.5mm} n \textrm{ odd}
\end{array}
\right., \qquad
\Phi_{n}^{ks}:=\phi_{n}^{ks}
+ B_0^{k}A_n^{ks},
$
\State
$
\begin{array}{l}
\mathcal{Q}_{0}^{ks}:=0\\
\\
\mathcal{Q}_{1}^{ks}:=\phi_{1}^{ks}
\end{array}, \quad
\mathcal{Q}_{2n}^{ks}:=
\begin{bsmallmatrix}
\phi_{2n}^{ks} \\ 
\overline{\Phi}_{2n-1}^{ks} \\ 
\Phi_{2n-2}^{ks} \\ 
\vdots  \\ 
\overline{\Phi}_{1}^{ks} 
\end{bsmallmatrix}, 
\quad
\mathcal{Q}_{2n+1}^{ks}:=
\begin{bsmallmatrix}
\phi_{2n+1}^{ks} \\ 
\overline{\Phi}_{2n}^{ks} \\ 
\overline{\Phi}_{2n-1}^{ks} \\ 
\vdots  \\ 
\Phi_{1}^{ks} 
\end{bsmallmatrix}, \quad
\psi_{n}^{krs}:=
\left\{\begin{array}{ll}
\mathcal{A}_n^{kr}\mathcal{Q}_n^{ks}, & b_{kr} \textrm{ odd}\\
\overline{\mathcal{A}}_n^{kr}\mathcal{Q}_n^{ks}, & b_{kr} \textrm{ even}
\end{array}
\right.,
$
\State 
$\xi_n^{krs}:=
\psi_n^{krs}+
\left\{
\begin{array}{ll}
\hspace{-2mm} (\overline{A}_n^{kr})^{T}B_0^{k}\overline{A}_0^{ks}, & \hspace{-2mm} b_{kr}, n \textrm{ odd} \\
\hspace{-2mm}(A_n^{kr})^{T}B_0^{k}A_0^{ks}, & \hspace{-2mm}  b_{kr}\textrm{ odd},n \geq 2 \textrm{ even}\\
\hspace{-2mm}(\overline{A}_n^{kr})^{T}B_0^{k}A_0^{ks},  & \hspace{-2mm} b_{kr}, n\geq 2 \textrm{ even}\\
\hspace{-2mm}(A_n^{kr})^{T}B_0^{k}\overline{A}_0^{ks},        &\hspace{-2mm}  b_{kr} \textrm{ even},n \textrm{ odd} \\
\hspace{-2mm}0, &  n=0
\end{array}
\right.$ \hspace{-2mm}
$
\Psi_{n}^{krs}:=\hspace{-0.5mm}\xi_{n}^{krs}
+\left\{\begin{array}{ll}
\hspace{-2mm}(A_0^{kr})^{T}B_0^{k}A_n^{ks}, & \hspace{-2mm} b_{kr} \textrm{ odd}\\
\hspace{-2mm}(\overline{A}_0^{kr})^{T}B_0^{k}A_n^{ks}, &\hspace{-2mm}  b_{kr} \textrm{ even}
\end{array}
\right.
$
\vspace{1mm}
\hspace{5mm}
\For {$j=0:\alpha_1-1$}
    \If {$r\in \{1,\ldots,N\}$,  $ j\in \{1,\ldots,\alpha_r-1\}$}
    \State $A_j^{rr}=\frac{1}{2}A_0^{rr}-\frac{1}{2}A_0^{rr}(C_0^{r})^{-1}\big(Z_j^{r}+\psi_j^{rrr}
    +\sum_{k=1}^{r-1}\Psi_{j-\alpha_k+\alpha_r}^{krr}+\sum_{k=r+1}^{N}\Psi_{j-\alpha_{r}+\alpha_k}^{krr}\big)$
    \EndIf
    \For {$p=1:N-1$}
        \If {$r\in \{1,\ldots,N\}$, $j\leq \alpha_{r+p}-1$, $r+p\leq N$}
        \small
        \State \hspace{-3mm} $A_j^{r(r+p)}=-A_0^{r(r+p)}(C_0^{r})^{-1}\big(
        \xi_j^{rr(r+p)}+\sum_{k=1}^{r-1}\Psi_{j-\alpha_k+\alpha_r}^{kr(r+p)}+ \sum_{k=r+1}^{r+p}\Psi_{j}^{kr(r+p)} $\\
        \qquad \qquad \qquad \qquad \qquad \qquad \qquad \qquad \qquad \qquad \qquad \qquad $+\sum_{k=r+p+1}^{N}\Psi_{j-\alpha_{r+p}+\alpha_k}^{kr(r+p)}\big)$   
\normalsize
        \EndIf
    \EndFor
\EndFor
\end{algorithmic}
For simplicity, in this algorithm we define $\sum_{j=l}^{n}a_j=0$ if $l>n$, and it is understood that the inner loop (i.e. for p =1 : N-1) is not performed for $N=1$. 


\end{lemma}

\vspace{-3mm}

To prove Lemma \ref{EqTca} we follow the same general approach as in \cite[Lemma 4.1]{TSOS},
however, some additional intrigueging technical problems arise.

For the sake of clarity we point out 
the correct order of calculating the entries of $\mathcal{X}$ in the lemmas. 
First, all nonzero entries of the blocks below the main diagonal of $\mathcal{X}=[\mathcal{X}_{rs}]_{r,s=1}^{N}$ (i.e. $A_j^{rs}$ for $r>s$) can be chosen freely.
We proceed by computing the upper triangular part of $\mathcal{X}$. We begin with the diagonal entries $A_0^{rr}$ of 
the main diagonal blocks 
$\mathcal{X}_{rr}$.
Next, the step $j=0$, $p=1$ (if $N\geq 2$) of the algorithm yields the diagonal entries of the first upper off-diagonal blocks of $\mathcal{X}$ (i.e. $(\mathcal{X}_{r(r+1)})_{11}=A_{0}^{r(r+1)}$). Further, the step $j=0$, $p=2$ gives the diagonal entries of the second upper off-diagonal blocks of $\mathcal{X}$ (i.e. $(\mathcal{X}_{r(r+2)})_{11}=A_{0}^{r(r+2)}$), 
and so forth. 
In the same fashion 
the step for fixed $j\in \{1,\ldots,\alpha_1-1\}$, $p\in \{0,\ldots,N\}$ yields the entries on the $j$-th upper off-diagonals of the $p$-th upper off-diagonal blocks of $\mathcal{X}$, i.e. $(\mathcal{X}_{r(r+p)})_{1(j+1)}=A_{j+1}^{r(r+p)}$ with $r+p\leq N$, $j\leq \alpha_{r+p}-1$.

\begin{proof}[Proof of Lemma \ref{EqT} (\ref{EqTII}) (\ref{EqTIII})]
We analyze (\ref{f1}) for $\mathcal{Y}=\mathcal{B}\mathcal{X}$, $\mathcal{\widetilde{X}}=\mathcal{F}\mathcal{X}^{T}\mathcal{F}$ (see (\ref{BBF}))
and $\mathcal{X}=[\mathcal{X}_{rs}]_{r,s=1}^{N}$ with $\mathcal{X}_{rs}$ as in (\ref{EqTX}). Observe that the fact
%
%
\begin{align}\label{ETE}
&E_{\alpha}(I_{n})\big(T(A_0,\ldots,A_{\alpha-1})\big)^{T}E_{\alpha}(I_{m})=T(A_0^{T},\ldots,A^{T}_{\alpha-1}), \quad A_0,\ldots,A_{\alpha-1}\in \mathbb{C}^{m\times n},
\end{align}
implies
\vspace{-2mm}
\begin{align*}
\widetilde{\mathcal{X}}_{rk}
=
E_{\alpha_r}(I_{m_r})\mathcal{X}_{kr}^{T} E_{\alpha_k}(I_{m_k})
&=
\left\{
\begin{array}{cc}
\begin{bsmallmatrix}
\widetilde{\mathcal{T}}_{rk}\\
0
\end{bsmallmatrix}, 
& \alpha_r >\alpha_k \\
\begin{bsmallmatrix}
0 & \widetilde{\mathcal{T}}_{rk}
\end{bsmallmatrix}, & \alpha_r<\alpha_k\\
\widetilde{\mathcal{T}}_{rk}, & \alpha_r=\alpha_k
\end{array}
\right., \,\,\,
\widetilde{T}_{rk}=
T\big((A_0^{kr})^{T},\ldots,(A_{b_{kr}-1}^{kr})^{T}\big).
\end{align*}

\vspace{-1mm}
For simplicity we set
$\Phi_{n}^{ks}:=\sum_{i=0}^{n} B_{n-i}^{k}A_i^{ks}$, $n\in \{0,\ldots, b_{rs}-1\}$, and we have
\begin{align}\label{YSP}
&\mathcal{Y}_{ks} 
=
\left\{
\begin{array}{ll}
\begin{bsmallmatrix}
\mathcal{S}_{ks}\\
0
\end{bsmallmatrix}, 
& \alpha_k >\alpha_s \\
\begin{bsmallmatrix}
0 & \mathcal{S}_{ks}
\end{bsmallmatrix}, & \alpha_k<\alpha_s\\
\mathcal{S}_{ks}, & \alpha_k=\alpha_s
\end{array}
\right., \quad  
\begin{array}{rl}
\mathcal{S}_{ks} &\hspace{-3mm} =T\big(B_0^{k},\ldots, B_{b_{ks}-1}^{k}\big)T\big(A_0^{ks},\ldots,A_{b_{ks}-1}^{ks}\big)\\
                 &\hspace{-3mm} =T\big(\Phi_0^{ks},\ldots, \Phi_{b_{ks}-1}^{ks}\big)
                 \end{array}.
\end{align}
Next, for $k,r,s\in \{1,\ldots,N\}$, $n\in \{0,\ldots, b_{rs}-1\}$ we set:
\begin{align}
\label{Prsk}
\Psi^{krs}_{n}:=
\left\{
\begin{array}{ll}
\hspace{-1mm} 
\begin{bsmallmatrix}
(A_0^{kr})^{T} & (A_1^{kr})^{T} & \ldots & (A_{n}^{rr})^{T} 
\end{bsmallmatrix}
\begin{bsmallmatrix}
\Phi_{n}^{ks} \\
\vdots \\
\Phi_{0}^{ks} \\
\end{bsmallmatrix}, & \hspace{-1mm} n\geq 0
\\
\hspace{-1mm} 0, & \hspace{-1mm} n<0
\end{array}
\right.
=
\left\{
\begin{array}{ll}
\hspace{-1mm}  \sum_{i=0}^{n}(A_i^{kr})^{T}\Phi_{n-i}^{ks}, & \hspace{-1mm}n\geq 0\\
\hspace{-1mm} 0, & \hspace{-1mm} n<0
\end{array}
\right.\hspace{-2mm}, 
\end{align}
\vspace{-3mm}
\begin{align}\label{simetrija}
(\Psi^{krs}_{n})^{T}
&=\sum_{i=0}^{n}(\Phi_i^{ks})^{T}A_{n-i}^{kr} 
=\sum_{i=0}^{n}\sum_{l=0}^{i}(A_l^{ks})^{T}( B_{i-l}^{k})^{T}A_{n-i}^{kr}
=\sum_{l=0}^{n}\sum_{i=l}^{n}(A_l^{ks})^{T} B_{i-l}^{k}A_{n-i}^{kr}\nonumber\\
&=\sum_{l=0}^{n}(A_{l}^{ks})^{T} \sum_{i'=0}^{n-l} B_{i'}^{k}A_{n-l-i'}^{kr}
=\sum_{l=0}^{n}(A_l^{ks})^{T}\Phi_{n-l}^{kr}=\Psi^{ksr}_{n}, \qquad n\geq 0.
\end{align}
%
Furthermore,
\vspace{-2mm}
\begin{align}\label{psijr}
(\widetilde{\mathcal{X}}_{rk})_{(1)}(\mathcal{Y}_{k(r+p)})^{(n+1)}
&=
\left\{
\begin{array}{ll}
\Psi_{n-\alpha_{r+p}+\alpha_k}^{kr(r+p)}, &  k\geq r+p+1\\
\Psi_{n}^{kr(r+p)},                       &  r+p \geq k\geq r+1, p\geq 1\\
\Psi_{n-\alpha_k+\alpha_r}^{kr(r+p)},     &  k\leq r\\
\end{array}
\right..
\end{align}

We now calculate matrices $A_0^{rr}$ for $r\in \{1,\ldots,N\}$. Since
\vspace{-1mm}
\[
(\widetilde{\mathcal{X}}_{rk})_{(1)}=\left\{
\begin{array}{ll}
\begin{bsmallmatrix}
(A_0^{kr})^{T} & * & \ldots & *
\end{bsmallmatrix}, & k\geq r\\
\begin{bsmallmatrix}
0 & * & \ldots & *
\end{bsmallmatrix}, & k<r
\end{array}
\right.,\qquad
(\mathcal{Y}_{kr})^{(1)}=
\left\{
\begin{array}{ll}
\begin{bsmallmatrix}
B_0^{k}A_0^{kr} \\
0\\
\vdots\\
0
\end{bsmallmatrix}, & k\leq r\\
0, & k>r
\end{array}
\right.,
\]
\vspace{-1mm}
we deduce $\sum_{k=1}^{N}(\widetilde{\mathcal{X}}_{rk})_{(1)}((\mathcal{Y})_{kr})^{(1)}=
(A_0^{rr})^{T}B_0^{r}A_0^{rr}$, thus 
(\ref{f1})
for $r=s$, $j=1$ yields:
\begin{equation}\label{GABA}
C_0^{r}=(A_0^{rr})^{T}B_0^{r}A_0^{rr}, \qquad r\in \{1,\ldots,N\}.
\end{equation}
If $B_0^{r},C_0^{r}$ are as in (\ref{EqTII}) (\ref{EqTIII}) then $\sqrt{\frac{v_0}{u_0}}I_{m_r}$ is one solution of (\ref{GABA}) of the form (\ref{AVW}).

\vspace{-1mm}
Proceed to the key step: an inductive computation of
the remaining entries.
Fix 
$p\in \{0,\ldots, N-1 \}$, $j\leq \alpha_r-1$, but not $p=j=0$. For $r,s,n$ satisfying
%
\begin{align}\label{induA}
&j\geq 1, n\in \{0,\ldots,j-1\}, s\geq r 
\quad\textrm{or}\quad p\geq 1, n=j, r\leq s\leq r+p-1 \\
&\textrm{ or }\quad  s\leq r, n\in \{0,\ldots,b_{rs}-1\}, N\geq 2, \nonumber
\end{align}
we assume that there exist $V_n^{rs}, W_n^{rs}\in \mathbb{C}^{m_r\times m_s}$ ($W_{-1}^{rs}:=0$, hence $F_{0}^{rs} =0$) so that
%
\vspace{-1mm}
\[
A_{n}^{rs}=\widetilde{A}_{n}^{rs}+F_{n}^{rs}, \qquad
\widetilde{A}_n^{rs}:=
\begin{bmatrix}
V_n^{rs} &  W_n^{rs} \\
-\mu^{2}\overline{W}_n^{rs} &  \overline{V}_n^{rs} 
\end{bmatrix},\quad
F_n^{rs}:=
\begin{bmatrix}
0 &  0 \\
\overline{W}_{n-1}^{rs} & 0 
\end{bmatrix}.
\]
We need to prove that $\widetilde{A}_j^{r(r+p)}=A_j^{r(r+p)}-F_{j}^{r(r+p)}$ is of the form 
$\begin{bsmallmatrix}
V &  W \\
-\mu^{2}\overline{W} &  \overline{V} 
\end{bsmallmatrix}$
as well. 

The trick of the proof is to reduce 
$(\mathcal{C}_{r(r+p)})_{1j}=((\widetilde{\mathcal{X}}\mathcal{Y})_{r(r+p)})_{1j}$ 
to a certain linear matrix equation in $\widetilde{A}_j^{r(r+p)}$ (and possi\-bly  $(\widetilde{A}_j^{r(r+p)})^{T}$) with coefficients of the appropriate form and depending only on $A_{n}^{rs}$ for $r,s,n$ satisfying (\ref{induA}).

If $n,r,s$ satisfy (\ref{induA}) or if $n=j$, $s=r+p$, we have ($L_{r}F_j^{rs}=0$, $(F_{j-1}^{kr})^{T}L_r=0$):
\vspace{-2mm}
\begin{align}\label{PhiDE}
\Phi_n^{rs} & =\sum_{i=0}^{n} B_{n-i}^{r}A_i^{rs}=K_{r} \sum_{i=0}^{n}u_{n-i}^{r}(\widetilde{A}_i^{rs}+F_i^{rs})+L_{r} \sum_{i=0}^{n-1}u_{n-1-i}^{r}(\widetilde{A}_i^{rs}+F_i^{rs})
\nonumber\\
         &= K_{r} D_n^{rs}+K_{r} E_{n}^{rs}+L_{r} D_{n-1}^{rs},\\
%
%
%
&D_{-1}^{rs}:=0, \quad D_n^{rs}:= 
\sum_{i=0}^{n}u_{n-i}^{r}\widetilde{A}_i^{rs},
\quad
E_n^{rs}:= 
\sum_{i=0}^{n}u_{n-i}^{r}F_i^{rs}.\nonumber
\end{align}
Further, we set
\vspace{-1mm}
\[
U_n^{rs}:=\hspace{-0.5mm}\sum_{i=0}^{n}u_{n-i}^{r}V_i^{rs}, \,\,\, Z_n^{rs}:=\hspace{-0.5mm}\sum_{i=0}^{n}u_{n-i}^{r}W_i^{rs}, \quad 
\big(
D_n^{rs}=\begin{bmatrix}U_n^{rs} & \hspace{-1mm}Z_n^{rs}\\
-\mu^{2}\overline{Z}_n^{rs} & \hspace{-1mm} \overline{U}_n^{rs}
\end{bmatrix},
E_n^{rs}=\begin{bmatrix}0 & \hspace{-1mm} 0\\
\overline{Z}_{n-1}^{rs} & \hspace{-1mm} 0
\end{bmatrix}
\hspace{-0.5mm}\big).
\]
Using this and (\ref{PhiDE}) it is straightforward to compute
\vspace{-1mm}
\begin{align*}
\Psi^{krs}_n  & =\sum_{i=0}^{n} (A_{i}^{kr})^{T}\Phi_{n-i}^{ks}
=  \sum_{i=0}^{n}
(\widetilde{A}_i^{kr})^{T}K_r D_{n-i}^{ks}
+\sum_{i=0}^{n-1}(\widetilde{A}_i^{kr})^{T} L_r D_{n-i-1}^{ks} 
+ \hspace{-1mm}\sum_{i=0}^{n-1}(\widetilde{A}_i^{kr})^{T}K_r E_{n-i}^{ks}\\
&\qquad \qquad \qquad\qquad \quad+ \hspace{-1mm}\sum_{i=1}^{n}(F_{i}^{kr})^{T}K_r D_{n-i}^{ks}
+ \hspace{-1mm}\sum_{i=1}^{n}(F_{i-1}^{kr})^{T}K_r E_{n-i}^{ks}=\\
%
=&
\sum_{i=0}^{n}
\begin{bsmallmatrix}
-\mu^{2}((V_{i}^{rs})^{T}U_{n-i} -\mu^{2}(\overline{W}_{i}^{rs})^{T}\overline{Z}_{n-i}) & \hspace{1mm}-\mu^{2}((V_{i}^{rs})^{T}Z_{n-i} +(\overline{W}_{i}^{rs})^{T}\overline{U}_{n-i})\\
-\mu^{2}((\overline{V}_{i}^{rs})^{T}\overline{Z}_{n-i} +(W_{i}^{rs})^{T}U_{n-i})   & (\overline{V}_{i}^{rs})^{T}\overline{U}_{n-i} -\mu^{2}(W_{i}^{rs})^{T}Z_{n-i}
\end{bsmallmatrix}\\
&+\sum_{i=0}^{n-1}\begin{bsmallmatrix}
-\mu^{2}(\overline{W}_{i}^{rs})^{T}\overline{Z}_{n-1-i}^{ks} & \quad (\overline{W}_{i}^{rs})^{T}\overline{U}_{n-1-i}^{ks}+(V_{i}^{rs})^{T}Z_{n-1-i}^{ks} \\
(\overline{V}_{i}^{rs})^{T}\overline{Z}_{n-1-i}+(W_{i}^{rs})^{T}U_{n-1-i}         & \quad (W_{i}^{rs})^{T}Z_{n-1-i}
\end{bsmallmatrix}\\
& +
\sum_{i=0}^{n-1}\begin{bsmallmatrix}
(V_{i}^{rs})^{T}U_{n-1-i} -\mu^{2}(\overline{W}_{i}^{rs})^{T}\overline{Z}_{n-1-i}^{ks} & 0 \\
0        & 0
\end{bsmallmatrix}
+
\sum_{i=0}^{n-2}\begin{bsmallmatrix}
(\overline{W}_{i}^{rs})^{T}\overline{Z}_{n-2-i}^{ks} & 0 \\
0        & 0
\end{bsmallmatrix}.
\end{align*}
%
Finally, we define
\begin{align}\label{gamakrs}
\Gamma_{-1}^{krs}:=0,\,\,\,\, 
\Gamma_n^{krs}:=&
\sum_{i=0}^{n}
\begin{bsmallmatrix}
-\mu^{2}((V_{i}^{rs})^{T}U_{n-i} -\mu^{2}(\overline{W}_{i}^{rs})^{T}\overline{Z}_{n-i}) & \hspace{1mm}-\mu^{2}((V_{i}^{rs})^{T}Z_{n-i} +(\overline{W}_{i}^{rs})^{T}\overline{U}_{n-i})\\
-\mu^{2}((\overline{V}_{i}^{rs})^{T}\overline{Z}_{n-i} +(W_{i}^{rs})^{T}U_{n-i})   & (\overline{V}_{i}^{rs})^{T}\overline{U}_{n-i} -\mu^{2}(W_{i}^{rs})^{T}Z_{n-i}
\end{bsmallmatrix}\\
&+\sum_{i=0}^{n-1}\begin{bsmallmatrix}
-\mu^{2}(\overline{W}_{i}^{rs})^{T}\overline{Z}_{n-1-i}^{ks} & \quad (\overline{W}_{i}^{rs})^{T}\overline{U}_{n-1-i}^{ks}+(V_{i}^{rs})^{T}Z_{n-1-i}^{ks} \\
(\overline{V}_{i}^{rs})^{T}\overline{Z}_{n-1-i}+(W_{i}^{rs})^{T}U_{n-1-i}         & \quad (W_{i}^{rs})^{T}Z_{n-1-i}
\end{bsmallmatrix}.\nonumber
\end{align}
Therefore, for $r,s,n$ satisfying (\ref{induA}) or for $n=j$, $s=r+p$ we can write
\begin{align}\label{PGD}
&
\Psi_{n}^{krs}=\Gamma_{n}^{krs}+
\begin{bsmallmatrix}
-\frac{1}{\mu^{2}}[\Gamma_{n-1}^{krs}]_{11} &  0\\
0                       &  0
\end{bsmallmatrix}.
%
\end{align}
%
%
%
%
%
Next, by applying (\ref{PGD}) and (\ref{psijr}) we further write;
$\Gamma_{n}^{kr(r+p)}:=0$ for $n<0$.):
\begin{align}
\sum_{k=1}^{N}(\widetilde{\mathcal{X}}_{rk})_{(1)}(\mathcal{Y}_{k(r+p)})^{(n+1)}
&=
\sum_{k=1}^{r}\Psi_{n-\alpha_k+\alpha_r}^{kr(r+p)}+
\sum_{k=r+1}^{r+p}\Psi_{n}^{kr(r+p)}
+\sum_{k=r+p+1}^{N}\Psi_{n-\alpha_{r+p}+\alpha_k}^{kr(r+p)}\nonumber\\ 
&=\gamma (n,r,p)+
\begin{bsmallmatrix}
-\frac{1}{\mu^{2}}[\gamma (n-1,r,p)]_{11} &  0\\
0                       &  0
\end{bsmallmatrix},
\label{vsotaX}
\end{align}
%
%
\vspace{-5mm}
\begin{align}\label{kapa1} 
\gamma (n,r,p) 
&:=
\Gamma_{n}^{kr(r+p)}+\left(\sum_{k=1}^{r-1}\Gamma_{n-\alpha_k+\alpha_r}^{kr(r+p)}+
\sum_{k=r+1}^{r+p}\Gamma_{j}^{kr(r+p)}
+\sum_{k=r+p+1}^{N}\Gamma_{n-\alpha_{r+p}+\alpha_k}^{kr(r+p)}\right).
\end{align}

Using (\ref{vsotaX}), the equation $((\widetilde{\mathcal{X}}\mathcal{Y})_{r(r+p)})_{1(j+1)}=(\mathcal{C}_{r(r+p)})_{1(j+1)}$ can be seen as 
\begin{equation*}
\gamma (j,r,p)
+\begin{bsmallmatrix}
-\frac{1}{\mu^{2}}[\gamma (j-1,r,p)]_{11} &  0\\
0                       &  0
\end{bsmallmatrix}
=\left\{\begin{array}{ll} v_{j}K_r+v_{j-1}L_r, & p=0\\ 0, & p\neq 0 \end{array}\right..
\end{equation*}
%
We show by induction that it is actually reduces to
%
\begin{equation}\label{induK}
\gamma (j,r,p)=\left\{\begin{array}{ll} v_{j}K_r, & p=0\\ 0, & p\neq 0 \end{array}\right..
\end{equation}
%
Indeed, it is clear for $j=0$ (since $v_{-1}=\gamma (-1,r,p)=0$), while 
assuming
(\ref{induK}) for some $n<j$ we easily conclude 
the following fact 
yielding the claim for $n+1$:
\[
\begin{bsmallmatrix}
[-\frac{1}{\mu^{2}}\gamma (n,r,p)]_{11} & 0\\
0                      & 0\\
\end{bsmallmatrix}
=\left\{
\begin{array}{ll} 
\begin{bsmallmatrix}
v_n[-\frac{1}{\mu^{2}} K_r]_{11} & 0\\
0                      & 0\\
\end{bsmallmatrix}, & p=0\\ 
0, & p\neq 0 \end{array}\right.
=
\left\{\begin{array}{ll} v_{n}L_r, & p=0\\ 0, & p\neq 0 \end{array}\right..
\]
%
Observe that $\Gamma^{rr(r+p)}_j$ (see (\ref{gamakrs})), $u_0^{r}(\widetilde{A}_0^{rr})^{T}K_r \widetilde{A}_{j}^{r(r+p)}$, $u_0^{r}(\widetilde{A}_j^{rr})^{T}K_r \widetilde{A}_{0}^{r(r+p)}$, 
and hen\-ce the expressions below
are both of the form 
$\begin{bsmallmatrix}
-\mu^{2}V &  W\\
\overline{W} & \overline{V}
\end{bsmallmatrix}$:
%
%
\begin{align*}
\Gamma_j^{rrr}-u_0^{r}(\widetilde{A}_0^{kr})^{T}K_r \widetilde{A}_{j}^{ks}+u_0^{r}(\widetilde{A}_j^{kr})^{T}K_r \widetilde{A}_{0}^{ks},\qquad 
\Gamma_j^{rr(r+p)}-u_0^{r}(\widetilde{A}_0^{kr})^{T}K_r \widetilde{A}_{j}^{ks}
.
\end{align*}
Moreover, the equation of 
(\ref{induK}) can be seen as:
%
\begin{align}
\big(u_0^{r}(A_0^{rr})^{T}K_r\big) \widetilde{A}_j^{r(r+p)}                                       &=\kappa(j,r,p), \qquad p\geq 1,\\
\big(u_0^{r}(A_0^{rr})^{T}K_{r}\big)\widetilde{A}_j^{rr}+(\widetilde{A}_j^{rr})^{T}\big(u_0^{r}K_{r}A_0^{rr}\big) &=\kappa(j,r,0), \qquad p=0,\nonumber
\end{align}
with $\kappa(j,r,p)$ of the form
$\begin{bsmallmatrix}
-\mu^{2}V &  W\\
\overline{W} & \overline{V}
\end{bsmallmatrix}$
and 
depending on $A_{n}^{rs}$ for $n,r,s$ satisfying (\ref{induA}).
Similarly, as we proved (\ref{simetrija}), 
we see that $\kappa(j,r,0)$ is symmetric.

Thus $\widetilde{A}_j^{r(r+p)}$ for $p\geq 1$ and $p=0$ is a solution of equations $A^{T}Y=B$ and $A^{T}X+X^{T}A=B$ for
$A=u_0^{r}(A_0^{rr})^{T}K_{r}$, $B=\kappa(j,r,p)$, respectively (i.e. $Y=(A^{T})^{-1}B$ and $X=(A^{T})^{-1}(\frac{1}{2}B+Z)$ with $Z$ skew-symmetric).
Since $(A^{T})^{-1}=((A_0^{rr})^{T}B_0^{r})^{-1}=
A_0^{r}(C_0^{r})^{-1}=\frac{1}{v_0^{r}}A_0^{r}(K_0^{r})^{-1}=
\frac{1}{v_0^{r}}
\begin{bsmallmatrix}
-\frac{1}{\mu^{2}}V_0^{rr} &  W_0^{rr}\\
\overline{W}_0^{rr} & \overline{V}_0^{rr}
\end{bsmallmatrix}$ (see (\ref{GABA})) and $B=\kappa(j,r,p)$ is of the form
$\begin{bsmallmatrix}
-\mu^{2}V &  W\\
\overline{W} & \overline{V}
\end{bsmallmatrix}$,
it follows that $Y$
is of the form  
$\begin{bsmallmatrix}
V &  W\\
-\mu^{2}\overline{W} & \overline{V}
\end{bsmallmatrix}$,
while $X$ is of this form
precisely when $Z$ is of this form. 
This completes the inductive step.
\end{proof}


\begin{proof}[Proof of Lemma \ref{EqTca}]
Let $\mathcal{X}=[\mathcal{X}_{rs}]_{r,s=1}^{N}$ with $\mathcal{X}_{rs}$ as in (\ref{EqTXca}) and $\mathcal{Y}=\mathcal{B}\mathcal{X}$, $\mathcal{\widetilde{X}}=\mathcal{F}\mathcal{X}^{T}\mathcal{F}$ (see (\ref{BBF})). Next, for $A_0,A_1,\ldots,A_{\alpha-1}\in \mathbb{C}^{m\times n}$ we have
%
\[
E_{\alpha}(I_{n})\big(T_c(A_0,A_1,\ldots,A_{\alpha-1})\big)^{T}E_{\alpha}(I_{m})
=\left\{
\begin{array}{ll}
T_c(\overline{A}_0^{T},A_1^{T},\ldots,
\overline{A}_{\alpha-2}^T,A_{\alpha-1}^T), & \alpha \textrm{ even}\\
T_c(A_0^{T},\overline{A}_1^{T},\ldots,
\overline{A}_{\alpha-2}^{T},A_{\alpha-1}^T), & \alpha \textrm{ odd}
\end{array}\right.;
\]
the entry in the first row and in the $j$-th column of the matrix $T_c(\overline{A}_0^{T},A_1^{T},\overline{A}_{2}^T,\ldots)$ (or $T_c(A_0^{T},\overline{A}_1^{T},\overline{A}_{2}^{T},\ldots)$) is $A_{j-1}$ for $j$ odd (even) and $\overline{A}_{j-1}$ for $j$ even (odd).
Thus  
\begin{align*}
&\widetilde{\mathcal{X}}_{rk}:=
E_{\alpha_r}(I_{m_r})\mathcal{X}_{sr}^{T} E_{\alpha_s}(I_{m_s})=
\left\{
\begin{array}{cc}
\begin{bmatrix}
\widetilde{\mathcal{T}}_{rk}\\
0
\end{bmatrix}, 
& \alpha_r >\alpha_k \\
\begin{bmatrix}
0 & \widetilde{\mathcal{T}}_{rk}
\end{bmatrix}, & \alpha_r<\alpha_k\\
\widetilde{\mathcal{T}}_{rk}, & \alpha_r=\alpha_k
\end{array}
\right.,
\end{align*}
\begin{align*}
&\widetilde{\mathcal{T}}_{rk}=
\left\{
\begin{array}{ll}
T_c\big((\overline{A}_0^{kr})^{T},(A_1^{kr})^{T},\ldots,(\overline{A}_{b_{kr}-2}^{kr})^{T},(A_{b_{kr}-1}^{kr})^{T}\big), & b_{kr} \textrm{ even}\\
T_c\big((A_0^{kr})^{T},(\overline{A}_1^{kr})^{T},\ldots,(\overline{A}_{b_{kr}-2}^{kr})^{T},(A_{b_{kr}-1}^{kr})^{T}\big), & b_{kr} \textrm{ odd}
\end{array}\right..
\end{align*}
We also have
\begin{align*}
&\mathcal{Y}_{ks} = 
\left\{
\begin{array}{cc}
\begin{bmatrix}
S_{ks}\\
0
\end{bmatrix}, 
& \alpha_k >\alpha_s \\
\begin{bmatrix}
0 & S_{ks}
\end{bmatrix}, & \alpha_k <\alpha_s\\
S_{ks}, & \alpha_k=\alpha_s
\end{array}
\right., \,\,\,
\begin{array}{rl}
S_{ks} &\hspace{-3mm} =T\big(B_0^{k},B_1^{k},\ldots, B_{b_{ks}-1}^{k}\big)T_c\big(A_0^{ks},A_1^{ks},\ldots,A_{b_{ks}-1}^{ks}\big)\\
       &\hspace{-3mm} =T_c\big(\Phi_0^{ks},\Phi_1^{ks},\ldots, \Phi_{b_{ks}-1}^{rs}\big)
\end{array}\hspace{-2mm},
\\
&\Phi_{2n}^{ks}\hspace{-1mm}
:=
\sum_{j=0}^{n}B_{2n-2j}^{k}A_{2j}^{ks}+\hspace{-1mm}\sum_{j=0}^{n-1}B_{2n-2j-1}^{k}\overline{A}_{2j+1}^{ks}, \,\,
\Phi_{2n+1}^{ks}\hspace{-1mm}
:=
\sum_{j=0}^{n}(B_{2n-2j}^{k}A_{2j+1}^{ks}+
B_{2n-2j+1}^{k}\overline{A}_{2j}^{ks}).
\end{align*}

Let us now compute $A^{rr}_0$ for $r\in \{1,\ldots,N \}$. Since
\vspace{-1mm}
\[
(\widetilde{X}_{rk})_{(1)}=\left\{
\begin{array}{ll}
\begin{bsmallmatrix}
(A_0^{rr})^{T} & * & \ldots *
\end{bsmallmatrix}, & k\geq r, \alpha_r \textrm{ odd}\\
\begin{bsmallmatrix}
(A_0^{rr})^{*} & * & \ldots * 
\end{bsmallmatrix}, & k\geq r, \alpha_r \textrm{ even}\\
\begin{bsmallmatrix}
0 & * & \ldots & *
\end{bsmallmatrix}, & k<r
\end{array}
\right.\hspace{-2mm},\quad
((\mathcal{Y})_{kr})^{(1)}=
\left\{
\begin{array}{ll}
\begin{bsmallmatrix}
B_0^{k}A_0^{kr} \\
0\\
\vdots\\
0
\end{bsmallmatrix}, & k\leq r\\
0, & k>r
\end{array}
\right.,
\]
\normalsize 
it follows from (\ref{f1}) for $r=s $, $j=1$ that
\vspace{-1mm}
\begin{equation}\label{GABAca}
C_0^{r}=
\left\{
\begin{array}{ll}
(A_0^{rr})^{T}B_0^{r}A_0^{rr}, & \alpha_r \textrm{ odd}\\
(A_0^{rr})^{*}B_0^{r}A_0^{rr}, & \alpha_r \textrm{ even}
\end{array}
\right..
\end{equation}
Since $B_0^{r},C_0^{r}$ are real symmetric, then by Sylvester's inertia theorem this equation for $\alpha_r$ even has a solution $A_0^{rr}$ precisely when $B_0^{r},C_0^{r}$ are of the same inertia.

Next, if $N\geq 2$, we fix arbitrarily the blocks below the main diagonal of $[\mathcal{X}_{rs}]_{r,s=1}^{N}$,
and then
inductively compute the remaining entries (as in the proof of Lemma \ref{EqT}).
We fix $p\in \{0,\ldots, N-1 \}$ and $j\leq \alpha_r-1$, but not $p=j=0$. 
To get $A_j^{r(r+p)}$ (step $j,p$ of the algorithm in (\ref{EqT3ca})), we solve $(\mathcal{C}_{r(r+p)})_{1j}=((\widetilde{\mathcal{X}}\mathcal{Y})_{r(r+p)})_{1j}$, while assuming that we have already determined matrices 
$A_{n}^{rs}$ for 
%
\begin{align}\label{induA1}
&j\geq 1, n\in \{0,\ldots,j-1\}, s\geq r 
\quad\textrm{or}\quad p\geq 1, n=j, r\leq s\leq r+p-1 \\
&\textrm{ or } \quad s\leq r, n\in \{0,\ldots,b_{rs}-1\}, N\geq 2, \qquad (1 \leq r, s\leq N).\nonumber
\end{align}

To simplify calculations we use $\mathcal{A}_0^{kr}$, $\Phi_{n}^{krs}$, $\Psi_{n}^{krs}$ defined in the algorithm in (\ref{EqT3ca}), and in addition we introduce the matrix vectors $\mathcal{P}_{n}^{ks}$ with $\Phi_{n-j+1}^{ks}$ (and $\overline{\Phi}_{n-j+1}^{ks}$) in the $j$-th row for $j\geq 2$ odd (even):
\vspace{-1mm}
\[
\Psi_{n}^{krs}=
\left\{\begin{array}{ll}
\mathcal{A}_n^{rk}\mathcal{P}_n^{ks} & b_{kr} \textrm{ odd}\\
\overline{\mathcal{A}}_n^{rk}\mathcal{P}_n^{ks} & b_{kr} \textrm{ even}
\end{array}
\right., \qquad
\mathcal{P}_{2n}^{ks}:=
\begin{bsmallmatrix}
\Phi_{2n}^{ks} \\ \overline{\Phi}_{2n-1}^{ks} \\ \vdots  \\ \overline{\Phi}_{1}^{ks} \\ \Phi_{0}^{ks} 
\end{bsmallmatrix}, 
%
\quad
\mathcal{P}_{2n+1}^{ks}:=
\begin{bsmallmatrix}
\Phi_{2n+1}^{ks} \\ \overline{\Phi}_{2n}^{ks} \\ \vdots  \\ \Phi_{1}^{ks} \\ \overline{\Phi}_{0}^{ks}
\end{bsmallmatrix},\quad n\geq 0.
\]
Further, for $n\geq 0$ we obtain:
\vspace{-2mm}
\small
\begin{align*}
(\mathcal{A}_{2n+1}^{kr}\overline{\mathcal{P}}_{2n+1}^{ks})^{T}
=&\sum_{j=0}^{n}(\overline{\Phi}_{2j+1}^{kr})^{T}A_{2n-2j}^{ks}+\sum_{j=0}^{n}(\Phi_{2j}^{kr})^{T}\overline{A}_{2n+1-2j}^{ks}
\\[-3pt]
=&\sum_{j=0}^{n}\sum_{l=0}^{j}\big((\overline{A}_{2l+1}^{kr})^{T}B_{2j-2l}^{k}+(A_{2l}^{kr})^{T}B_{2j+1-2l}^{k}\big)A_{2n-2j}^{ks}\\[-1pt]
&+\big(\sum_{j=0}^{n}\sum_{l=0}^{j}(A_{2l}^{kr})^{T}B_{2j-2l}^{k}+\sum_{j=1}^{n}\sum_{l=0}^{j-1}(\overline{A}_{2l+1}^{kr})^{T}B_{2j-1-2l}^{k}\big)\overline{A}_{2n+1-2j}^{ks}
\\[-3pt]
=& \sum_{l=0}^{n}\sum_{j=l}^{n}(\overline{A}_{2l+1}^{kr})^{T}B_{2j-2l}^{k}A_{2n-2j}^{ks}+\sum_{l=0}^{n}\sum_{j=l}^{n}(A_{2l}^{kr})^{T}B_{2j+1-2l}^{k}A_{2n-2j}^{ks}\\[-3pt]
&+\sum_{l=0}^{n}\sum_{j=l}^{n}(A_{2l}^{kr})^{T}B_{2j-2l}^{k}\overline{A}_{2n+1-2j}^{ks}+\sum_{l=0}^{n-1}\sum_{j=l+1}^{n}\vspace{-1mm}(\overline{A}_{2l+1}^{kr})^{T}B_{2j-1-2l}^{k}\overline{A}_{2n+1-2j}^{ks}\\[-3pt]
=&(\overline{A}_{2n+1}^{kr})^{T}B_0^{k}A_0^{ks}
+\sum_{l=0}^{n}(A_{2l}^{kr})^{T}\sum_{j'=0}^{n-l}\big(B_{2j'+1}^{k}A_{2n-2l-2j'}^{ks}+B_{2j'}^{k}\overline{A}_{2n+1-2l-2j'}^{ks}\big)\\[-3pt]
&+\sum_{l=0}^{n-1}(\overline{A}_{2l+1}^{kr})^{T}\big(B_0^{k}A_{2n-2l}^{ks}+\hspace{-1mm}\sum_{j'=0}^{n-1-l}\hspace{-1mm}(B_{2j'+2}^{k}A_{\scriptscriptstyle{2n-2l-2j'-2}}^{ks} +B_{2j'+1}^{k}\overline{A}_{\scriptscriptstyle{2n-1-2l-2j}}^{ks})\big)
\\[-3pt]
=&\sum_{l=0}^{n}(\overline{A}_{2l+1}^{kr})^{T} \Phi_{2n-2l}^{ks}+ \sum_{l=0}^{n}(A_{2l}^{kr})^{T} \overline{\Phi}_{2n+1-2l}^{ks}=\mathcal{A}_{2n+1}^{kr}\overline{\mathcal{P}}_{2n+1}^{ks},
\end{align*}
%
\small
\vspace{-2mm}
\begin{align*}
(\overline{\mathcal{A}}_{2n}^{kr}\mathcal{P}_{2n}^{ks})^{T}
=&\sum_{j=1}^{n}(\overline{\Phi}_{2j-1}^{kr})^{T}A_{2n+1-2j}^{ks}+\sum_{j=0}^{n}(\Phi_{2j}^{kr})^{T}\overline{A}_{2n-2j}^{ks}\\[-1pt]
=&\sum_{j=1}^{n}\sum_{l=0}^{j-1}\big((A_{2l}^{kr})^{T}B_{2j-1-2l}^{k}+(\overline{A}_{2l+1}^{kr})^{T}B_{2j-2-2l}^{k}\big)A_{2n+1-2j}^{ks}\\[-1pt]
&+\big(\sum_{j=0}^{n}\sum_{l=0}^{j}(A_{2l}^{kr})^{T}B_{2j-2l}^{k}+\sum_{j=1}^{n}\sum_{l=0}^{j-1}(\overline{A}_{2l+1}^{kr})^{T}B_{2j-1-2l}^{k}\big)\overline{A}_{2n-2j}^{ks}\\[-3pt]
=&\sum_{l=0}^{n-1}\sum_{j=l+1}^{n}\big((A_{2l}^{kr})^{T}B_{2j-1-2l}^{k}A_{2n+1-2j}^{ks}\big)+\sum_{l=0}^{n}\sum_{j=l}^{n}\big((A_{2l}^{kr})^{T}B_{2j-2l}^{k}\overline{A}_{2n-2j}^{ks}\big)\\[-1pt]
&+\sum_{l=0}^{n-1}\sum_{j=l+1}^{n}\big((\overline{A}_{2l+1}^{kr})^{T}B_{2j-2l}^{k}A_{2n+1-2j}^{ks}+(\overline{A}_{2l+1}^{kr})^{T}B_{2j-1-2l}^{k}\overline{A}_{2n-2j}^{ks}\big)
\end{align*}
\begin{align*}
=&(A_{2n}^{kr})^{T}B_0^{k}\overline{A}_0^{ks}
+\sum_{l=0}^{n-1}(\overline{A}_{2l+1}^{kr})^{T}\sum_{j'=0}^{n-l-1}\big(B_{2j'-2}^{k}A_{2n-1-2j-2l}^{ks}+B_{2j'+1}^{k}\overline{A}_{2n-2j'-2l-2}^{ks}\big)\\[-3pt]
+&\sum_{l=0}^{n-1} (A_{2l}^{kr})^{T}\left(B_{0}^{k}\overline{A}_{2n-2l}^{ks}+\sum_{j'=0}^{n-l-1}\big(B_{2j'+1}^{k}A_{2n-1-2j'-2l}^{ks}+B_{2j'+2}^{k}\overline{A}_{2n-2j'-2l-2}^{ks}\big)\right)\\[-3pt]
=&\sum_{l=0}^{n}(A_{2l}^{kr})^{T} \overline{\Phi}_{2n-2l}^{ks}+ \sum_{l=0}^{n-1}(\overline{A}_{2l+1}^{kr})^{T} \Phi_{2n-1-2l}^{ks}=\mathcal{A}_{2n}^{kr}\overline{\mathcal{P}}_{2n}^{ks}.
\end{align*}
\normalsize
In a similar manner we prove
\vspace{-1mm}
\[
(\mathcal{A}_{2n+1}^{kr}\mathcal{P}_{2n+1}^{ks})^{T}=\overline{\mathcal{A}}_{2n+1}^{kr}\overline{\mathcal{P}}_{2n+1}^{ks}, \qquad (\mathcal{A}_{2n}^{kr}\mathcal{P}_{2n}^{ks})^{T}=\mathcal{A}_{2n}^{kr}\mathcal{P}_{2n}^{ks}.
\]
The above computations thus yield
\vspace{-1mm}
\begin{equation}\label{rssim}
(\Psi_{j}^{krs})^{T}=
\left\{\begin{array}{ll}
\Psi_j^{ksr} & j - b_{kr} \textrm{ odd}\\
\overline{\Psi}_j^{ksr} & j - b_{kr} \textrm{ even}
\end{array}
\right..
\end{equation}

\vspace{-2mm}
Since
\vspace{-1mm}
\begin{align*}
&(\widetilde{X}_{rr})_{(1)}=
\left\{
\begin{array}{ll}
\mathcal{A}_{\alpha_r-1}^{rr}, &  \alpha_r \textrm{ odd}\\
\overline{\mathcal{A}}_{\alpha_r-1}^{rr}, & \alpha_r \textrm{ even}
\end{array}
\right., 
\quad
(\mathcal{Y}_{r(r+p)})^{(j+1)}
=
\left\{\begin{array}{ll}
\mathcal{P}_{\alpha_r-1}^{rr}, & p=0,j=\alpha_r-1\\
\begin{bsmallmatrix}
\mathcal{P}_{j}^{r(r+p)}\\
0
\end{bsmallmatrix},   &  j< \alpha_{r}-1 \textrm{ or } p\geq 1
\end{array}
\right.\hspace{-1mm},
\end{align*}
%
we have
$(\widetilde{\mathcal{X}}_{rr})_{(1)}(\mathcal{Y}_{r(r+p)})^{(j+1)}=
\Psi_{j}^{rr(r+p)}.$
%
In particular, for $\xi(j,r,p)$ as defined in the algorithm in (\ref{EqT3ca}), we deduce for $j\geq 1$, $ p\geq 1$ and for $p=0$, respectively:
\vspace{-1mm}
\begin{align}\label{xijrp}
&(\widetilde{\mathcal{X}}_{rr})_{(1)}(\mathcal{Y}_{r(r+p)})^{(j+1)}
=
\Psi_{j}^{rr(r+p)}
=\xi(j,r,p)+\left\{\begin{array}{ll}
\hspace{-1mm}(A_0^{rr})^{T}B_0^{r}A_j^{r(r+p)}, & \hspace{-1mm}\alpha_{r} \textrm{ odd}\\
\hspace{-1mm}(A_0^{rr})^{*}B_0^{r}A_j^{r(r+p)}, & \hspace{-1mm}\alpha_{r} \textrm{ even},
\end{array}
\right.
\end{align}
\vspace{-3mm}
\begin{align*}
(\widetilde{\mathcal{X}}_{rr})_{(1)}&(\mathcal{Y})_{rr}^{(j+1)}
= 
\Psi_{j}^{rrr}
=\xi(j,r,0)+
\left\{\begin{array}{ll}
\hspace{-1mm}(\overline{A}_0^{rr})^{*}B_0^{r}A_j^{rr}+(A_j^{rr})^{*}B_0^{r}\overline{A}_0^{rr}, &  \hspace{-1mm}\alpha_r,j \textrm{ odd}\\
\hspace{-1mm}(A_0^{rr})^{T}B_0^{r}A_j^{rr}+ (A_j^{rr})^{T}B_0^{r}A_0^{rr}, &  \hspace{-1mm}\alpha_s\textrm{ odd},j \textrm{ even}\\
\hspace{-1mm}(A_0^{rr})^{*}B_0^{r}A_j^{rr}+(A_j^{rr})^{*}B_0^{r}A_0^{rr}, &  \hspace{-1mm}\alpha_r,j \textrm{ even}\\
\hspace{-1mm}(\overline{A}_0^{rr})^{T}B_0^{r}A_j^{rr}+(A_j^{rr})^{T}B_0^{r}\overline{A}_0^{rr}, &  \hspace{-1mm}\alpha_r\textrm{ even},j\textrm{ odd} \\
\end{array}
\right. \nonumber
\end{align*}

Summands of the second term in (\ref{f1})
for $j+1$ instead of $j$ consist
of
\vspace{-1mm}
\begin{align*}
&(\widetilde{X}_{rk})_{(1)}
=
\left\{
\begin{array}{ll}
\mathcal{A}_{\alpha_r-1}^{kr}, &  \alpha_r \textrm{ odd}\\
\overline{\mathcal{A}}_{\alpha_r-1}^{kr}, &  \alpha_r \textrm{ even}
\end{array}
\right., \quad 
%
(\mathcal{Y}_{k(r+p)})^{(j+1)}=
\left\{\begin{array}{ll}
\mathcal{P}_j^{r(r+p)}, & k=r+p\\
\begin{bsmallmatrix}
\mathcal{P}_j^{r(r+p)} \\
0
\end{bsmallmatrix}, & k< r+p ,\\
\begin{bsmallmatrix}
\mathcal{P}_{j-\alpha_{r+p}-\alpha_k}^{r(r+p)} \\
0
\end{bsmallmatrix}, & k> r+p
,
\end{array}
\right.,
\end{align*}
%
hence (for $N\geq r+1\geq 2$):
\vspace{-1mm}
\begin{align}\label{thjrp}
\Theta(j,r,p):=
& \sum_{k=r+1}^{N}(\widetilde{\mathcal{X}}_{rk})_{(1)}(\mathcal{Y}_{k(r+p)})^{(j+1)}\\
=
&\left\{\begin{array}{ll}
\sum_{k=r+1}^{N}\Psi_{j-\alpha_{r}+\alpha_k}^{krr}, & j\geq 1, p=0\\
\sum_{k=r+1}^{r+p}\Psi_{j}^{kr(r+p)}
+\sum_{k=r+p+1}^{N}\Psi_{j-\alpha_{r+p}+\alpha_k}^{kr(r+p)}, & j\geq 0, p\geq 1. 
\end{array}
\right.\nonumber
\end{align}
For simplicity, we defined $\sum_{k=r+p+1}^{N}\Psi_{j-\alpha_{r+p}-\alpha_k}^{rr(r+p)}=0$ for $r+p+1>N$.

Finally, the third term in (\ref{f1}) for $j+1$ instead of $j$ (with $N\geq 2$, $k\leq r-1$) is
\vspace{-1mm}
\begin{equation}\label{lajrp}
\Lambda(j,r,p):=\sum_{k=1}^{r-1}(\widetilde{\mathcal{X}}_{rk})_{(1)}(\mathcal{Y}_{k(r+p)})^{(j+1)}=\sum_{k=1}^{r-1}\Psi_{j-\alpha_k+\alpha_r}^{kr(r+p)},
\end{equation}
since we have
\[
(\widetilde{\mathcal{X}}_{rk})_{(1)}
=
\left\{
\begin{array}{ll}
\begin{bsmallmatrix}
0 &\mathcal{A}_{\alpha_{k}}^{kr} 
\end{bsmallmatrix}, & \alpha_{k} \textrm{ odd}\\
\begin{bsmallmatrix}
0 &\overline{\mathcal{A}}_{\alpha_{k}}^{kr} 
\end{bsmallmatrix}, & \alpha_{k} \textrm{ even}
\end{array}\right.,\qquad
(\mathcal{Y}_{k(r+p)})^{(j+1)}
=\begin{bsmallmatrix}
\mathcal{P}_j^{k(r+p)} \\
0
\end{bsmallmatrix}, \quad 1\leq k\leq r-1.
\]

For $j,p\geq 0$ with $j+p\geq 1$ we define
\begin{equation}\label{Djrp}
D_j^{r(r+p)}:=\Xi(j,r,p)+\Theta(j,r,p)+\Lambda(j,r,p).
\end{equation}
%
We combine $(\mathcal{C}_{r(r+p)})_{1j}=((\widetilde{\mathcal{X}}\mathcal{Y})_{r(r+p)})_{1j}$ in (\ref{f1}) with (\ref{xijrp}), (\ref{thjrp}), (\ref{lajrp}), (\ref{Djrp}):
\begin{align}
\label{eqATB}
&(A_0^{rr})^{*}B_0^{r}A_j^{r(r+p)}=-D_j^{r(r+p)}, \qquad
\alpha_{r} \textrm{ even}, \qquad p\geq 1\\
&(A_0^{rr})^{T}B_0^{r}A_j^{r(r+p)}=-D_j^{r(r+p)}, \qquad
\alpha_{r} \textrm{ odd}, \qquad p\geq 1,\nonumber
\\
\label{eqATTAB}
&(\overline{A}_0^{rr})^{*}B_0^{r}A_j^{rr}+(A_j^{rr})^{*}B_0^{r}\overline{A}_0^{rr}=C_j^{r}-D_j^{rr}, \qquad \alpha_r,j \textrm{ odd}, (p=0),j\geq 1\\
&(A_0^{rr})^{T}B_0^{r}A_j^{rr}+(A_j^{rr})^{T}B_0^{r}A_0^{rr}=C_j^{r}-D_j^{rr}, \qquad \alpha_r \textrm{ odd}, j \textrm{ even } (p=0),j\geq 1\nonumber\\
&(A_0^{rr})^{*}B_0^{r}A_j^{rr}+(A_j^{rr})^{*}B_0^{r}A_0^{rr}=C_j^{r}-D_j^{rr}, \qquad \alpha_r,j \textrm{ even }, (p=0),j\geq 1\nonumber\\
&(\overline{A}_0^{rr})^{T}B_0^{r}A_j^{rr}+(A_j^{rr})^{T}B_0^{r}\overline{A}_0^{rr}=C_j^{r}-D_j^{rr}, \qquad \alpha_r \textrm{ even }, j \textrm{ odd} (p=0),j\geq 1  \nonumber
\end{align}
Moreover, from (\ref{rssim}) it follows that $\Psi^{krs}_{n}$ for $r=s$, and thus $\xi(j,r,0)$, $\Theta(j,r,0)$, $\Lambda(j,r,0)$, $C_j^{r}-D_j^{rr}$ are all symmetric (Hermitian) if $\alpha_r-j$ is odd (even).

Since (\ref{GABAca}) is equivalent to 
$A_0^{r}(C_0^{r})^{-1}=\left\{
\begin{array}{ll}
((A_0^{rr})^{T}B_0^{r})^{-1}, & \alpha_{r} \textrm{ odd}\\
((A_0^{rr})^{*}B_0^{r})^{-1}, &  \alpha_{r} \textrm{ even}
\end{array}
\right.$, 
(\ref{eqATB}) yields $A_j^{r(r+p)}=-A_0^{r}(C_0^{r})^{-1}D_j^{r(r+p)}$ for $p\geq 1$.
Next, we get $A_j^{rr}$ by solving (\ref{eqATTAB}), i.e. an equation of the form $A^{T}X+X^{T}A=B$ for $\alpha_r-j$ odd and of the form $A^{*}X+X^{*}A=B$ for $\alpha_r-j$ even, with given $A$ nonsingular and $B$ symmetric or Hermitian; 
the solution in the first case is $X=\frac{1}{2}(A^{T})^{-1}B+(A^{T})^{-1}Z$ with $Z$ skew-symmetric and in the second case $X=\frac{1}{2}(A^{*})^{-1}B+(A^{*})^{-1}Z$ with $Z$ skew-Hermitian.
If $\alpha_r$ is odd (even), then for $j$ even (odd) we have $A=B_0^{r} A_0^{rr}$ ($A=B_0^{r} \overline{A}_0^{rr}$), hence
$(A^{T})^{-1}=
A_0^{r}(C_0^{r})^{-1}$,
while a similar argument for 
$\alpha_r-j$ even 
gives $(A^{*})^{-1}=A_0^{r}(C_0^{r})^{-1}$. 
Furthermore, 
$B=C_j^{r}-D_j^{rr}$ and it
depends only on $A_n^{rs}$ with $n,r,s$ satisfying (\ref{induA1}). It is straightforward to conclude the algorithm in (\ref{EqT3ca}).

It is only left to sum up the dimensions:
\vspace{-2mm}
\small
\begin{align*}
&2\sum_{r=1}^{N}\sum_{s=1}^{r-1}\alpha_s m_r m_s
+ \sum_{\alpha_r \textrm{ even}}\big(m_r^{2}+\tfrac{(\alpha_r-2)m_r^{2}}{2}+\tfrac{\alpha_r}{2}m_r(m_r-1)\big)\\
&+\sum_{\alpha_r \textrm{ odd}}\big(\tfrac{1}{2}m_r(m_r-1)+\tfrac{\alpha_r-1}{2}m_r^{2}+\tfrac{\alpha_r-1}{2}m_r(m_r-1)\big)\\
=&\sum_{r=1}^{N}\big(\alpha_r m_r^{2} +2\sum_{s=1}^{r-1}\alpha_s m_r m_s\big) - \sum_{\alpha_r \textrm{ even}}\tfrac{1}{2}m_r\alpha_r-\sum_{\alpha_r \textrm{ odd}}\tfrac{1}{2} m_r(\alpha_r+1).
\end{align*}
\normalsize
This completes the proof of the lemma.
\end{proof}

\begin{remark} \begin{enumerate}[label=\arabic*.,ref={\arabic*},wide=0pt,itemsep=2pt] 
\item
It would be interesting to find a nice description of $A_0^{rr}$ of the form (\ref{AVW}) and such that $C_0^{r }=(A_0^{rr})^{T}B_0^{r}A_0^{rr}$ with $B_0^{r }$, $C_0^{r }$ as in (\ref{aass}).
\item
One could consider (\ref{eqFYFIY}) even for nonsingular $\mathcal{B}$ and $\mathcal{C}$, since the solutions of $A^{T}X+X^{T}A=B$ and $A^{*}X+X^{*}A=B$ in this case are known (see \cite{Brad}, \cite{LR}).
\end{enumerate}
\end{remark}


\begin{example}
We solve (\ref{eqFYFIY}) for $\mathcal{F}=E_3(I)\oplus E_2(I)$, $\mathcal{B}=\mathcal{B}'=I_6(I)$ with the identity matrix $I$, and where the solution $\mathcal{X}_c$ is of the form as in Example \ref{exX}. 
We have:
\small
\begin{align*}
&\qquad
\widetilde{\mathcal{X}}_c\mathcal{X}_c=\begin{bmatrix}[ccc|cc]
A_1^{T} & B_1^{*} & C_1^{T}  &  N_1^{*}  &  P_1^{T} \\  
0   & A_1^{*} & B_1^{T}     &  0  &   N_1^{T} \\ 
0   & 0   & A_1^{T}      &  0 & 0   \\  
\hline
0    & H_1^{*} & F_1^{T}                        &  A_2^{*}  &  B_2^{T} \\  
0   & 0     & H_1^{T}              &  0    &   A_2 ^{T}  
\end{bmatrix}
\begin{bmatrix}[ccc|cc]
A_1 & B_1 & C_1  &  H_1  &  F_1 \\  
0   & \overline{A}_1 & \overline{B}_1     &  0  &   \overline{H}_1  \\  
0   & 0   & A_1     &  0 & 0   \\  
\hline
0     & N_1 & P_1                         &  A_2  &  B_2 \\  
0    & 0   & \overline{N}_1              &  0    &   \overline{A}_2    
\end{bmatrix}=\\
&=
\begin{bmatrix}[ccc|cc]
A_1^{T}A_1 & A_1^{T}B_1+B_1^{*}\overline{A}_1 & A_1^{T}C_1+C_1^{T}A_1  &  A_1^{T}H_1+N_1^{*}A_2  &  A_1^{T}F_1+B_1^{*}\overline{H}_1 \\  
           &   +  N_1^{*} N_1      & + B_1^{*}\overline{B}_1  +N_1^{*}P_1+P_1^{T}\overline{N}_1                    &                           &      +N_1^{*}B_2     +P_1^{T}\overline{A}_2   \\  
  0         & A_1^{*}\overline{A}_1 &  A_1^{*}\overline{B}_1+B_1^{T}A_1 +N_1^{T}\overline{N}_1    &  0  &   A_1^{*}\overline{H}_1+N_1^{T}\overline{A}_2 \\ 
0   &     0                  & A_1^{T}A_1                        &  0  & 0             \\  
\hline 
   &    &                     &  A_{2}^{*}A_2  &   A_2^{*}B_2+B_2^{T}\overline{A}_2 \\  
   &    &                       &                &    + H_1^{*}\overline{H}_1       \\  
   &    &             &   0      &   A_2^{T}\overline{A}_2       
\end{bmatrix}
\end{align*}
\normalsize
By comparing diagonals of the main diagonal blocks in 
$\widetilde{\mathcal{X}}_c\mathcal{X}_c=\mathcal{I}$, we deduce that $A_1$ is orthogonal, while $A_2$ is unitary. Next, we choose $N_1$, $P_1$ arbitra\-rily. The diagonal element of the right upper block gives $A_1^{*}H_1+N_1^{*}A_2=0$, thus $H_1=-\overline{A}_1^{}N_1^{*}A_2$. 

We observe the first upper diagonals of the blocks to get $A_1^{T}B_1+B_1^{*}\overline{A}_1+N_1^{*}N_1=0$, $A_2^{*}B_2+B_2^{T}\overline{A}_2+H_1^{*}\overline{H}_1=0$ and $A_1^{T}F_1+B_1^{*}\overline{H}_1+N_1^{*}B_2 +P_1^{T}\overline{A}_2=0$. Thus $B_1=-\frac{1}{2}A_1N^{*}_1N_1+A_1Z_1$, $B_2=-\frac{1}{2}A_2 H^{*}_1\overline{H}_1+A_2 Z_2=-\frac{1}{2}N_1 N^{T}_1\overline{A}_2+A_2 Z_2$ for any $Z_1=-Z_1^{*}$, $Z_2=-Z_2^{T}$, and $F_1=-A_1(B_1^{*}\overline{H}_1+N_1^{*}B_2 +P_1^{T}\overline{A}_2)$. Finally, the second upper diagonal of the left upper block yields $A_1^{T}C_1+C_1^{T}A_1+B_1^{*}B_1 +N_1^{*}P_1+P_1^{T}\overline{N}_1=0$, therefore $C_1$ follows.
%
\end{example}

\vspace{-1mm}

Solutions of (\ref{eqFYFIY}) with $\mathcal{C}=\mathcal{B}$ form a group. Indeed, for any pair of solutions $\mathcal{X}_1,\mathcal{X}_2$ the product $\mathcal{X}_1\mathcal{X}_2^{-1}$ is a solution, too:
%
\begin{align*}
\mathcal{F}(\mathcal{X}_1\mathcal{X}_2^{-1})^{T}\mathcal{F}\mathcal{B} (\mathcal{X}_1\mathcal{X}_2^{-1})
&=\mathcal{F}(\mathcal{X}_2^{-1})^{T}\mathcal{F}\mathcal{F}\mathcal{X}_1^{T}\mathcal{F}\mathcal{B} \mathcal{X}_1\mathcal{X}_2^{-1}
=\mathcal{F}(\mathcal{X}_2^{-1})^{T}\mathcal{F}\mathcal{B}\mathcal{X}_2^{-1}=\\
&=\mathcal{F}(\mathcal{X}_2^{-1})^{T}\mathcal{F}\mathcal{B}(\mathcal{B}^{-1}\mathcal{F}\mathcal{X}_2^{T}\mathcal{F}\mathcal{B})=\mathcal{B}.
\end{align*}
%
%
Generators of this group are relatively simple as described below.

\begin{lemma}\label{lemauni}
Assume $\mathbb{T}_{}^{\alpha,\mu}$ and $\mathbb{T}_{c}^{\alpha,\mu}$ with $\alpha=(\alpha_1,\ldots,\alpha_n)$, $\mu=(m_1,\ldots,m_N)$ are as in (\ref{0T0}).
Let $\mathbb{X}_{}\subset \mathbb{T}_{}^{\alpha,\mu}$ and $\mathbb{X}_{c}\subset \mathbb{T}_{c}^{\alpha,\mu}$ be the sets of solutions $[\mathcal{X}_{rs}]_{r=1}^{N}$ with $\mathcal{X}_{rs}$ of the form (\ref{EqTX}) and of the form (\ref{EqTXca}), respectively, 
of the equation (\ref{eqFYFIY}) for $\mathcal{C}=\mathcal{B}$.
Then 
\vspace{-1mm}
\[
\mathbb{X}_{}=\mathbb{O}_{}\ltimes \mathbb{V}_{}\subset \mathbb{T}_{}^{\alpha,\mu}, \qquad \mathbb{X}_{c}=\mathbb{O}_{c}\ltimes \mathbb{V}_{c}\subset \mathbb{T}_{c}^{\alpha,\mu},
\]
in which the group $\mathbb{O}_{}$ (the group $\mathbb{O}_{c}$) consists of all matrices
$\mathcal{Q}=\bigoplus_{r=1}^{N}\big(\bigoplus_{j=1}^{\alpha_r} Q_{r}\big)$ ($\mathcal{Q}=\bigoplus_{r=1}^{N}( Q_{r}\oplus \overline{Q}_r\oplus Q_r\oplus \cdots)$) for $Q_{r}\in \mathbb{C}^{m_r\times m_r}$ such that $B_0^{r}=Q_{r}^{T}B_0^{r}Q_{r}$ (such that $B_0^{r}=Q_{r}^{T}B_0^{r}Q_{r}$ for $\alpha_r$ odd and $B_0^{r}=Q_{r}^{*}B_0^{r}Q_{r}$ for $\alpha_r$ even), $B_0^{r}=[\mathcal{B}_{rr}]_{11}$, while any $\mathcal{V}\in \mathbb{V}_{}$ (any $\mathcal{V}\in \mathbb{V}_{c}$)
can be written as 
$\mathcal{V}=\prod_{j=0}^{n_{\mathcal{V}}}\mathcal{V}_j$, where 
$\mathcal{V}_0=\bigoplus_{r=1}^{N}\mathcal{W}_r$ with $\mathcal{W}_r$ (complex-alternating) upper unitriangular Toeplitz 
and $\mathcal{V}_1,\ldots,\mathcal{V}_n$ of the form 
(\ref{0T0}) 
with (\ref{Hptk}). Both, $\mathbb{V}_{}$ and $\mathbb{V}_{c}$, are unipotent of order at most $\leq \alpha_1-1$.
Furthermore:

\vspace{-1mm}
\begin{enumerate}
\item \label{lemauni1}
If $\mathcal{B}$ is of the form (\ref{aass}) and $\mathcal{V}\in \mathbb{V}$ is of the form
(\ref{EqTX}) with (\ref{AVW}),
then $\mathcal{V}_0,\mathcal{V}_1,\ldots,\mathcal{V}_n$ can be choosen of the form 
(\ref{EqTX}) with (\ref{AVW}) as well.

\item \label{lemauni2}
If $\mathcal{B}=\bigoplus_{r=1}^{N}\big(\bigoplus_{j=1}^{\alpha_r} B_0^{r}\big)$, then $\mathbb{V}_{}$
is generated by matrices of the form (\ref{asZ}) and of the form (\ref{0T0}) with (\ref{Hptk}), (\ref{propS}),
while $\mathbb{V}_{c}$
is generated by matrices of the form (\ref{asZ2}) and of the form (\ref{0T0}) with (\ref{Hptk}), (\ref{propS2}) for $B_r=B_0^{r}$. 
\end{enumerate}
\end{lemma}

\vspace{-1mm}
If solutions of (\ref{eqFYFIY}) consist of
rectangular upper triangular Toeplitz blocks, the lemma coincides with \cite[Lemma 4.2]{TSOS}. Its proof works mutatis mutandis for solutions 
with rectangular complex-alternating upper triangular Toeplitz blocks, 
and also
for 
$\mathcal{B}$ of the form (\ref{aass}) and $\mathcal{X}$ of the form (\ref{EqTX}) with (\ref{AVW}).

\vspace{-1mm}
\section{Proofs of Theorem \ref{stabw} and Theorem \ref{stabz}}\label{sec2}

To get the isotropy group at $\mathcal{H}^{\varepsilon}$ we shall find all orthogonal $Q$ that solve 
\vspace{-1mm}
\begin{equation}\label{HQQH}
\mathcal{H}^{\varepsilon}\overline{Q}= Q\mathcal{H}^{\varepsilon}.
\end{equation}
We shall first apply Lemma \ref{lemaBHH} to obtain a general solution of (\ref{HQQH}) (Proposition \ref{resAoXXA} (\ref{resAoXXA2})). 
It will then be written in a suitable form by using permutation ma\-tri\-ces from Lemma \ref{lemaP}.
Fi\-na\-lly, we take into account the orthogonality of solutions,
which yi\-elds to the crux of the problem, i.e. the 
equation (\ref{eqFYFIY}) considered in Sec. \ref{cereq}.
Applying Lemma \ref{EqT} and Lemma \ref{EqTca} will thus immediately imply Therem \ref{stabw}, while further using Lemma \ref{lemauni} will furnish Theorem \ref{stabz}.

\begin{enumerate}[label={\bf Case \arabic*.},ref={Case \arabic*},wide=0pt,itemsep=15pt]

\item[{\bf Case \ref{stabz01}}.] 
Suppose
\vspace{-2mm}
\[
\mathcal{H}^{\varepsilon}
=\bigoplus_{r=1}^{N}\big( \bigoplus_{j=1}^{m_r} \varepsilon_{r,j} H_{\alpha_r}(\lambda)\big),\qquad \rho:=\lambda^{2},\quad\lambda\geq 0,\quad \textrm{all }\varepsilon_{r,j}\in \{-1,1\},
\]
where $H_{\alpha_r}(\lambda)$ is as in (\ref{Hmz}) for $z=\lambda$, $m=\alpha_r$.
We have
\vspace{-2mm}
\begin{align}\label{HTS}
\mathcal{H}^{}:=
\bigoplus_{r=1}^{N}\big( \bigoplus_{j=1}^{m_r}  H_{\alpha_r}(\lambda)\big)=S_{\varepsilon}\mathcal{H}^{\varepsilon}\overline{S}_{\varepsilon}^{-1}, \qquad
S_{\varepsilon}=\bigoplus_{r=1}^{N}\big(\bigoplus_{j=1}^{m_r}\sqrt{\varepsilon_{r,j}}I_{\alpha_j}\big).  
\end{align}
Using (\ref{HTS}), the equation (\ref{HQQH}) further transforms to
%
\begin{equation}\label{H1QQH1}
\mathcal{H}^{}\overline{Y}= Y\mathcal{H}^{},\qquad Y=S_{\varepsilon}Q S_{\varepsilon}^{-1}.
\end{equation}
Lemma \ref{lemaBHH} (\ref{BHH2}) gives the solution $Y\hspace{-0.8mm}=\hspace{-0.8mm}P^{-1}XP$ of (\ref{H1QQH1}), so the solution of (\ref{HQQH}) is
\[
Q=S_{\varepsilon}^{-1}P^{-1}XPS_{\varepsilon}, \qquad
P=\bigoplus_{r=1}^{N}\big(\bigoplus_{j=1}^{m_r}P_{\alpha_r}\big), \quad
P_{\alpha}:=\tfrac{e^{-i\frac{\pi}{4}}}{\sqrt{2}}(I_{\alpha}+iE_{\alpha}),
\]
in which $X\hspace{-0.5mm}=\hspace{-0.5mm}[X_{rs}]_{r,s=1}^{N}$ is such that $X_{rs}$ is an $m_r$-by-$m_s$ block matrix with blocks 
of the form (\ref{QTY}) for $m=\alpha_r$, $n=\alpha_s$ and $T$ is an $b_{rs}$-by-$b_{rs}$ real (complex-alterna\-ting) upper triangular Toeplitz matrix for $\lambda>0$ ($\lambda=0$); $b_{rs}=\min\{\alpha_r,\alpha_s\}$.

\quad
Since $P_{\alpha}=P_{\alpha}^T
$, 
$P_{\alpha}^{2}=
E_{\alpha}$, 
we get 
$P^{2} =(P^{-1})^{2}=E:=\bigoplus_{r=1}^{N}\left(\bigoplus_{j=1}^{m_r}  E_{\alpha_r}  \right)$. Next, $ S_{\varepsilon} =S_{\varepsilon}^T$,
$S_{\varepsilon}^2=(S_{\varepsilon}^2)^{-1}=\bigoplus_{r=1}^{N}\big(\bigoplus_{j=1}^{m_r}\varepsilon_{r,j} I_{\alpha_j}\big)$, $PS_{\varepsilon}^2=S_{\varepsilon}^2P$.
Thus $I=Q^TQ$
becomes
%
\begin{align}\label{QTQH}
I=   & \big(S_{\varepsilon}^TP^TX^T(P^{-1})^T(S_{\varepsilon}^{-1})^T\big)\big(S_{\varepsilon}^{-1}P^{-1}X P S_{\varepsilon}\big)\nonumber\\
 I  = & PS_{\varepsilon} \big(S_{\varepsilon}^T P^TX^T(P^{-1})^T(S_{\varepsilon}^{-1})^T S_{\varepsilon}^{-1}P^{-1}X P S_{\varepsilon}\big) S_{\varepsilon}^{-1}P^{-1}\\
 I =   &  S_{\varepsilon}^2 P^2 X^T (P^{-1})^{2}S_{\varepsilon}^{-2}X\nonumber\\
    S_{\varepsilon}^{2} = & E X^{T} ES_{\varepsilon}^{2}X.\nonumber
\end{align}
%
%
We conjugate matrices of (\ref{QTQH}) by 
$\Omega=\bigoplus_{r=1}^{N}\Omega_{\alpha_r,m_r}$ 
from Lemma \ref{lemaP}:
%
\begin{align}\label{ortoC3}
 \Omega^{T}S_{\varepsilon}^{2}\Omega = & (\Omega^{T}E\Omega)(\Omega^{T} Y^{T} \Omega)(\Omega^{T}E\Omega)(\Omega^{T}S_{\varepsilon}^{2}\Omega )(\Omega^{T}Y\Omega)\\
 \mathcal{B} = &\mathcal{F}\mathcal{X}^{T}\mathcal{F}\mathcal{B}\mathcal{X},\nonumber
\end{align}
where $\mathcal{F}\hspace{-0.5mm}=\hspace{-0.5mm}\Omega^{T}E\Omega=\hspace{-0.5mm}\bigoplus_{r=1}^{N}E_{\alpha_r}(I_{m_r})$, \hspace{-0.5mm}$\mathcal{B}=\Omega^{T}S_{\varepsilon}^{2}\Omega=\hspace{-0.5mm}\bigoplus_{r=1}^{N}\hspace{-0.5mm}\big(\bigoplus_{k=1}^{\alpha_r}(\bigoplus_{j=1}^{m_r}\varepsilon_{r,j} )\big)$ and $\mathcal{X}\hspace{-0.5mm}=\hspace{-0.5mm}\Omega^{T}X\Omega$ for $\lambda>0$ (for $\lambda=0)$ is of the form
(\ref{0T02})
with real (complex-alter\-nating) upp\-er tri\-an\-gu\-lar Toeplitz blocks.
Lemma \ref{EqT} (\ref{EqT2}), (\ref{EqT3}), (\ref{EqTII}) (\ref{EqTIIi}), Lemma \ref{EqTca} and Lemma \ref{lemauni} (\ref{lemauni2}) give Theorem \ref{stabw} for $\rho\geq 0$ and Theorem \ref{stabz} (\ref{stabz01}).

%
%
%
\item[{\bf Case \ref{stabz2}}.] 

Let
\[
\mathcal{H}^{\varepsilon} 
=\bigoplus_{r=1}^{N}\big( \bigoplus_{j=1}^{m_r} K_{\alpha_r}(\mu)\big), \qquad \rho:=-\mu^{2},\quad \mu> 0,
\]
where $K_{\alpha_r}(\mu)$ is as in (\ref{KLmz}) for $z=\mu$, $m=\alpha_r$.
Lemma \ref{lemaBHH} (\ref{BHHK}) now solves (\ref{HQQH}): 
\begin{equation}\label{QK}
Q=P^{-1}V^{-1}SXS^{-1}VP,
\end{equation}
in which $X=[X_{rs}]_{r,s=1}^{N}$ with an $m_r$-by-$m_s$ block matrix $X_{rs}$ whose blocks are
of the form (\ref{QTC})
for
$T_{1}$, $T_{2}$ 
of the form (\ref{QTY}) for $m=\alpha_r$, $n=\alpha_s$ and $T$ upper triangular Toeplitz of size $b_{rs}\times b_{rs}$ with $b_{rs}=\min\{\alpha_r,\alpha_s\}$, and
%
\begin{align*}
&P=\bigoplus_{r=1}^{N}\big(\bigoplus_{j=1}^{m_r}e^{\frac{i\pi}{4}}(P_{\alpha_r}\oplus P_{\alpha_r})\big), \qquad
V=\bigoplus_{r=1}^{N} \big(\bigoplus_{j=1}^{m_r}e^{i\frac{\pi}{4}}(W_{\alpha_r}\oplus \overline{W}_{\alpha_r})\big), \nonumber
\\
& S=\bigoplus_{r=1}^{N}\big(\bigoplus_{k=1}^{m_r}\begin{bsmallmatrix}
0 & U_{\alpha_r} \\
J_{\alpha_r}(-i\mu)\overline{U}_{\alpha_r} & 0
\end{bsmallmatrix}\big),\qquad 
P_{\alpha }:=\tfrac{e^{-i\frac{\pi}{4}}}{\sqrt{2}}(I_{\alpha }+iE_{\alpha }), \quad W_{\alpha }:= \oplus_{j=0}^{\alpha -1} i^{j},\nonumber
\end{align*}
where $U_{\alpha }$ is a solution of
the equation 
$U_{\alpha }J_{\alpha }(-\mu^{2})=(J_{\alpha }(i\mu))^{2}U_{\alpha }$.
Observe that $P=P^T
$, $P^{2} =-(P^{-1})^{2}=iE$ with $E:=\bigoplus_{r=1}^{N}\big(\bigoplus_{j=1}^{m_r} (E_{\alpha_r}\oplus E_{\alpha_r}) \big)$ and $V=V^{T}$, $V^{-1}=\overline{V}$. 
If we define $B=-iES^{T}\overline{V}E\overline{V}S$, then $I=Q^TQ$ is equivalent to 
%
\small
\begin{align}\label{eqSTIS}
I=   & \big(P^{T}V^{T}(S^{-1})^{T}X^{T}S^{T}(V^{-1})^{T}(P^{T})^{-1}\big)\big(P^{-1}V^{-1}S X S^{-1}VP\big)\nonumber\\
(S^{T}V^{-1}P^{-1})( P^{-1}V^{-1}S) = & S^{T}V^{-1}P^{-1}  \big( PV (S^{-1})^{T}X^{T}S^{T}\overline{V}(-iE)\overline{V}S X S^{-1}VP \big) P^{-1}V^{-1}S\nonumber\\
 -iES^{T}\overline{V}E\overline{V}S = & EX^{T}E(-iES^{T}\overline{V}E\overline{V}S) X\nonumber\\
 B = & EX^{T}EB X.
\end{align}
\normalsize

\quad
Next, since $(J_{\alpha_r}(-i\mu))^{T}E_{\alpha_r}=E_{\alpha_r}J_{\alpha_r}(-i\mu)$,
we have
\vspace{-1mm}
\begin{align*}
\overline{U}_{\alpha_r}^{T}(J_{\alpha_r}(-i\mu))^{T}E_{\alpha_r}J_{\alpha_r}(-i\mu)\overline{U}_{\alpha_r}
&=
\overline{U}_{\alpha_r}^{T}E_{\alpha_r}(J_{\alpha_r}(-i\mu))^{2}\overline{U}_{\alpha_r}=
\overline{U}_{\alpha_r}^{T}E_{\alpha_r}\overline{U}_{\alpha_r}J_{\alpha_r}(-\mu^{2}).
\end{align*}
We combine it with a calculation
$\overline{V}E\overline{V}=\oplus_{r=1}^{N}\big( i^{\alpha_r}\oplus_{j=1}^{m_r} ( (-1)^{\alpha_r}E_{\alpha_r}\oplus -E_{\alpha_r})\big)$:
\vspace{-1mm}
\small
\begin{equation}\label{IU}
B=\bigoplus_{r=1}^{N}\Big( i^{\alpha_r-1}\bigoplus_{j=1}^{m_r} 
\begin{bsmallmatrix}
-E_{\alpha_r}\overline{U}_{\alpha_r}^{T}E_{\alpha_r}\overline{U}_{\alpha_r}J_{\alpha_r}(-\mu^{2}) & 0 \\
0  & (-1)^{\alpha_r}E_{\alpha_r}U_{\alpha_r}^{T}E_{\alpha_r}U_{\alpha_r}
\end{bsmallmatrix}\Big).
\end{equation}
\normalsize
Furthermore, we show that
$E_{\alpha_r}U_{\alpha_r}^{T}E_{\alpha_r}U_{\alpha_r}$ is upper triangular Toplitz:
\vspace{-1mm}
\begin{align*}
\big( E_{\alpha_r}U_{\alpha_r}^{T}E_{\alpha_r}U_{\alpha_r}\big)J_{\alpha_r}(-\mu^{2})=&E_{\alpha_r}U_{\alpha_r}^{T}E_{\alpha_r}(J_{\alpha_r}(i\mu))^{2}U_{\alpha_r}\\
=&E_{\alpha_r}U_{\alpha_r}^{T}((J_{\alpha_r}(i\mu))^{2})^{T}E_{\alpha_r}U_{\alpha_r}\\
=&E_{\alpha_r}(U_{\alpha_r}J_{\alpha_r}(-\mu^{2}))^{T}E_{\alpha_r}U_{\alpha_r}\\
=&J_{\alpha_r}(-\mu^{2})\big(E_{\alpha_r}U_{\alpha_r}^{T}E_{\alpha_r}U_{\alpha_r}\big).
\end{align*}
We choose $U_{\alpha_r}$ so that
the odd (even) rows have real (purely imaginary or zero) entries, e.g. 
$U_{\alpha_r}=
\begin{bmatrix}
v_1 & v_2 & \ldots & v_{\alpha_r}
\end{bmatrix}$ is formed by taking real eigenvector $v_0$ of $(J_{\alpha_r}(i\mu))^{2}$ and then recursively solve equations $((J_{\alpha_r}(i\mu))^{2}+\mu^{2})v_{n}=v_{n-1}$ for $n\in \{2,\ldots,\alpha_r\}$. All nonvanishing entries of $E_{\alpha_r}U_{\alpha_r}^{T}E_{\alpha_r}U_{\alpha_r}$ are hence purely imaginary for $\alpha_r$ even and real for $\alpha_r$ odd.
Up to real scaling $U_{\alpha_r}$, we deduce
\vspace{-1mm}
\begin{equation*}
B=\bigoplus_{r=1}^{N}\Big(\bigoplus_{j=1}^{m_r} 
\begin{bsmallmatrix}
T(1,u_1^{r},\ldots,u_{\alpha_r-1}^{r})J_{\alpha_r}(-\mu^{2}) & 0 \\
0  & T(1,u_1^{r},\ldots,u_{\alpha_r-1}^{r})
\end{bsmallmatrix}\Big), \qquad  u_1^{r},\ldots,u_{\alpha_r-1}^{r}\in \mathbb{R}.
\end{equation*}

\vspace{-2mm}
\quad
Proceed by conjugating (\ref{eqSTIS}) with $\Omega_{}'=\bigoplus_{r=1}^{N}\Omega_{\alpha_r,m_r}'$ 
as in Lemma \ref{lemaP} (\ref{lemaP2}): 
\vspace{-1mm}
\begin{align}\label{ortoD2}
 (\Omega')^{T} B\Omega' = & \big((\Omega')^{T}E\Omega'\big)\big((\Omega')^{T} X\Omega'\big)^{T}\big((\Omega')^{T}E\Omega'\big)\big((\Omega')^{T}B \Omega'\big)\big((\Omega')^{T} X\Omega'\big)\\
  \mathcal{B} = &\mathcal{F}'\mathcal{X}^{T}\mathcal{F}'\mathcal{B}\mathcal{X}; \nonumber
\end{align}
where $\mathcal{F}'\hspace{-0.6mm}=(\Omega')^{T}E\Omega'\hspace{-0.6mm}=\hspace{-0.6mm}\bigoplus_{r=1}^{N}\hspace{-1mm}E_{\alpha_r}(I_{2m_r})$, \hspace{-0.6mm}$
\mathcal{B}\hspace{-0.6mm}=\hspace{-0.6mm}(\Omega')^{T}B\Omega'\hspace{-0.6mm}=\hspace{-0.6mm}\bigoplus_{r=1}^{N}\hspace{-1mm}T(B_0^{r},\ldots,B_{\alpha_r-1}^{r}) 
$
with $B_n^{r}$ as in (\ref{aass}), 
and $\mathcal{X}=(\Omega')^{T}X\Omega'$ of the form (\ref{0T0}) with (\ref{Trsvw0}) for $\rho=-\mu^{2}$.
To prove Theorem \ref{stabz} (\ref{stabz2}), we apply Lemma \ref{EqT} (\ref{EqT2}), (\ref{EqT3}), (\ref{EqTII}) (\ref{EqTIII}) and Lemma \ref{lemauni} (\ref{lemauni1}) to (\ref{ortoD2}), while to conclude the proof of Theorem \ref{stabw} for $\rho<0$,
it remains to find $\dim (\Sigma_{\mathcal{H}})$,
since Lemma \ref{EqT} does not provide it in this case.

\quad
We directly compute the dimension of the tangent space of $\Orb(\mathcal{H}^{\varepsilon})$
which is diffeomorphic to the quotient of the orthogonal group over $\Sigma_{\mathcal{H}^{\varepsilon}}$ (\cite[Ch. II.1]{GOV}).
If $Q(t)$ is a complex-differentiable path of orthogonal matrices with $Q(0)=I$, then differentiation of $(Q(t))^{T}Q(t)=I$ at $t=0$ yields $Z_0:=\frac{d}{dt}\big|_{t=0}Q(t)=-\frac{d}{dt}\big|_{t=0}Q^{T}(t)=-Z_0^{T}$, and the tangent vector of the orbit at $\mathcal{H}^{\varepsilon}$ is
\vspace{-1mm}
\[
\frac{d}{dt}\big|_{t=0}\big(Q^{*}(t)\mathcal{H}^{\varepsilon} Q(t)\big)
=\tfrac{d}{dt}\big|_{t=0}Q^{*}(t)\mathcal{H}^{\varepsilon}+\mathcal{H}^{\varepsilon} \tfrac{d}{dt}\big|_{t=0}Q(t)
=-\overline{Z}_0\mathcal{H}^{\varepsilon}+\mathcal{H}^{\varepsilon}Z_0;
\]
$e^{tZ}$ is orthogonal for $Z=-Z^{T}$ 
with 
$\frac{d}{dt}\big|_{t=0}e^{tZ}=Z$.
Hence the codimension of $\Sigma_{\mathcal{H}^{\varepsilon}}$ in the set of orthogonal matrices is equal to the codimension of $\{-\overline{Z}\mathcal{H}^{\varepsilon}+\mathcal{H}^{\varepsilon}Z=0\hspace{-1mm}\mid Z=-Z^{T}\}$ in the space of skew-symmetric matrices.
We must thus find those $Q$ in (\ref{QK}) (solving (\ref{HQQH}))
that satisfy $Q=-Q^{T}$. 
By recalling (\ref{eqSTIS}) with $P=P^{T}$, $P^{-2}=E$, $V=V^{T}$, $V^{-1}=\overline{V}$ and $B=-iES^{T}\overline{V}E\overline{V}S$, we deduce:
%
\begin{align}\label{eqSTIS2}
P^{-1}V^{-1}SXS^{-1}VP=  &  -P^{T}V^{T}(S^{T})^{-1}X^{T}S^{T}(V^{-1})^{T}(P^{-1})^{T}\nonumber\\
ES^{T}(V^{T})^{-1}(P^{T})^{-1}P^{-1}V^{-1}SX
=  & 
-
EX^{T}S^{T}(V^{-1})^{T}(P^{-1})^{T}P^{-1}V^{-1}S\\
BX= &-EX^{T}EB.\nonumber
\end{align} 
In the same manner as we transformed (\ref{eqSTIS}) to (\ref{ortoD2}), we transform (\ref{eqSTIS2}) to
%
\begin{align}\label{BXXB}
\mathcal{B}\mathcal{X} & =-\mathcal{F}'\mathcal{X}^{T}\mathcal{F}'\mathcal{B} \nonumber\\
\qquad \qquad \qquad \mathcal{B}_{rr}\mathcal{X}_{rs} & =-E_{\alpha_r}(I_{2m_r})\mathcal{X}_{sr}^{T}E_{\alpha_s}(I_{2m_s})\mathcal{B}_{ss}, \qquad r,s\in \{1,\ldots,N\}.
\end{align}
Clearly, $\mathcal{X}_{sr}$ for $r\neq s$ is uniquely determined by $\mathcal{X}_{rs}$. 
We now examine the case $r=s$.
We compare the entries in the first row of the $(j+1)$-th column
in (\ref{BXXB}):
%
\vspace{-1mm}
\begin{equation}\label{sBAAB}
\sum_{n=0}^{j}B_n^{r}A_{j-n}^{rr} =-\sum_{n=0}^{j}(A_{j-n}^{rr})^{T}B_n^{r}, \qquad r \in \{1,\ldots,N\}.
\end{equation}

\vspace{-2mm}
Since 
$B_n^{r }=(-\mu^{2}u_n^{r}+u_{n-1}^{r})I_{m_r}\oplus u_n^{r}I_{m_r}$,
$A_n^{rr}=
\begin{bsmallmatrix}
V_{n}^{rr} & W_{n}^{rr}\\
-\mu^{2}\overline{W}_{n}^{rr}+\overline{W}_{n-1}^{rr} & \overline{V}_n^{rr}
\end{bsmallmatrix}
$ with
$V_{n}^{rr},W_{n}^{rr}\in \mathbb{C}^{m_r\times m_r}$ for $n\in\{0,\ldots,\alpha_r-1\}$ and $u_{-1}^{r}=0$, $W_{-1}^{rr}=0$ (see (\ref{aass}), (\ref{AVW})),
then (\ref{sBAAB}) for $j=0$ gives
$
\begin{bsmallmatrix}
-\mu^{2}V_{0}^{rr} & -\mu^{2}W_{0}^{rr}\\
-\mu^{2}\overline{W}_{0}^{rr} & \overline{V}_0^{rr}
\end{bsmallmatrix}
=-
\begin{bsmallmatrix}
-\mu^{2}(V_{0}^{rr})^{T} & -\mu^{2}(\overline{W}_{0}^{rr})^{T}\\
-\mu^{2}(W_{0}^{rr})^{T} & (\overline{V}_0^{rr})^{T}
\end{bsmallmatrix}
$,
while for $j\geq 1$ it yields: 
\vspace{-1mm}
\begin{align}\label{BAABV}
\sum_{n=0}^{j}u_n^{r} 
\begin{bsmallmatrix}
-\mu^{2}V_{j-n}^{rr} & -\mu^{2}W_{j-n}^{rr}\\
-\mu^{2}\overline{W}_{j-n}^{rr} & \overline{V}_{j-n}^{rr}
\end{bsmallmatrix}  
+
\sum_{n=0}^{j-1}
u_n^{r}
\begin{bsmallmatrix}
 V_{j-1-n}^{rr} &  W_{j-1-n}^{rr}\\
 \overline{W}_{j-1-n}^{rr} & 0
\end{bsmallmatrix}
& = \\
=-\sum_{n=0}^{j}
u_n^{r}
\begin{bsmallmatrix}
-\mu^{2}(V_{j-n}^{rr})^{T} & -\mu^{2}(\overline{W}_{j-n}^{rr})^{T}\\
-\mu^{2}(W_{j-n}^{rr})^{T} & (\overline{V}_{j-n}^{rr})^{T}
\end{bsmallmatrix} 
& -
\sum_{n=0}^{j-1}
u_n^{r}
\begin{bsmallmatrix}
 (V_{j-1-n}^{rr})^{T} &   (W_{j-1-n}^{rr})^{T}\\
 (\overline{W}_{j-1-n}^{rr})^{T} & 0
\end{bsmallmatrix}. \nonumber
\end{align}
We prove by induction that 
$V_j^{rr}=-(V_j^{rr})^{T}$, $W_j^{rr}=-(W^{rr}_j)^{*}$ for all $j$. 
Clearly, $V_0^{rr}=-(V_0^{rr})^{T}$, $W_0^{rr}=-(W^{rr}_0)^{*}$. If we assume that the statement holds for  
$n < j$, it then follows 
from (\ref{BAABV}) that
$
\begin{bsmallmatrix}
-\mu^{2}V_{j}^{rr} & -\mu^{2}W_{j}^{rr}\\
-\mu^{2}\overline{W}_{j}^{rr} & \overline{V}_{j}^{rr}
\end{bsmallmatrix}
=
-
\begin{bsmallmatrix}
-\mu^{2}(V_{n-1}^{rr})^{T} & -\mu^{2}(\overline{W}_{n-1}^{rr})^{T}\\
-\mu^{2}(W_{n-1}^{rr})^{T} & (\overline{V}_{n-1}^{rr})^{T}
\end{bsmallmatrix}
$,
thus $V_{j}^{rr}=-(V_{j}^{rr})^{T}$, $W_{j}^{rr}=-(W^{rr}_{j})^{*}$. It remains to count all free parameters.


\item[{\bf Case \ref{stabz3}}.] 

Let
\[
\mathcal{H}^{\varepsilon} 
=\bigoplus_{r=1}^{N}\Big( \bigoplus_{j=1}^{m_r}  L_{\alpha_r}(\xi)\Big), \qquad \rho:=\xi^{2}\in \mathbb{C}\setminus \mathbb{R};
\]
$L_{\alpha_r}(\xi)$ is as in (\ref{KLmz}) for $z=\xi$, $m=\alpha_r$.
Lemma \ref{lemaBHH} (\ref{BHH2}) gives the solution of (\ref{HQQH}): 
%
\begin{equation}\label{caseLQ}
Q=P^{-1}XP,\qquad
P=\bigoplus_{r=1}^{N}\big(\bigoplus_{j=1}^{m_r}P_{\alpha_r}\oplus P_{\alpha_r}\big), \quad
P_{\alpha}:=\tfrac{e^{-i\frac{\pi}{4}}}{\sqrt{2}}(I_{\alpha}+iE_{\alpha}),
\end{equation}
where
$X=[X_{rs}]_{r,s=1}^{N}$ such that $X_{rs}$ is an $m_r$-by-$m_s$ block matrix whose blocks are 
of the form (\ref{QTC})
for $m=\alpha_r$, $T_{2}=0$ and 
$T_{1}$
of the form (\ref{QTY}) for $m=\alpha_r$, $n=\alpha_s$ with $T$ an $b_{rs}$-by-$b_{rs}$ complex upper triangular Toeplitz; $b_{rs}=\min\{\alpha_r,\alpha_s\}$.

\quad 
Similarly, (\ref{QTQH}) was obtained, we now apply (\ref{caseLQ}) to $I=Q^TQ$ to deduce 
%
\begin{align*}
I=   & P^TX^T(P^{-1})^TP^{-1}X P \\
    I = & E X^{T} EX,\nonumber
\end{align*}
in which $E:=\bigoplus_{r=1}^{N}\left(\bigoplus_{j=1}^{m_r} (E_{\alpha_r}\oplus E_{\alpha_r}) \right)$.
Using $\Omega_0$ from Lemma \ref{lemaP} (\ref{lemaP2}) 
we get
\begin{align}\label{ortoD3}
 I = & (\Omega_0^{T}E\Omega_0)(\Omega_0^{T} X\Omega_0)^{T}(\Omega_0^{T}E\Omega_0)(\Omega_0^{T} X\Omega_0)\nonumber\\
 I= & (\mathcal{F}\oplus  \mathcal{F}) \mathcal{X}^{T}  (\mathcal{F}\oplus  \mathcal{F})\mathcal{X},\\
 I= &\mathcal{F}\mathcal{V}^{T}\mathcal{F}\mathcal{V}, \nonumber
\end{align}
%
in which $\mathcal{F}=\bigoplus_{r=1}^{N}E_{\alpha_r}(I_{m_r})$, $\mathcal{X}=\Omega_0^{T} X\Omega_0=\mathcal{V}\oplus \overline{\mathcal{V}}$ for $\mathcal{V}$ of the form (\ref{0T0}) with upper tri\-angu\-lar Toeplitz blocks.
Finally, we apply Lemma \ref{EqT} (\ref{EqT2}), (\ref{EqT3}) and Lemma \ref{lemauni} (\ref{lemauni2})
to prove Theorem \ref{stabw} for $\rho\in \mathbb{C}\setminus \mathbb{R}$ and Theorem \ref{stabz} (\ref{stabz2}).
\end{enumerate}

\quad
This concludes the proof of the theorems.

\begin{remark}\label{conclR}
\begin{enumerate}[label=\arabic*.,ref={\arabic*},wide=0pt,itemsep=2pt] 
\item \hspace{-2mm} Solvability of (\ref{HQQH}) 
was first studied by the author \cite[Eq. 2.12]{TSH} to prove the uniqueness of Hong's normal form under orthogonal *congruence. The technique used there
was developed in \cite[Lemma 4.1]{TSOS} to the extent of solving (\ref{ortoD3}), and finally in this paper we give a complete solution of (\ref{HQQH}).
\item 
\hspace{-2mm} By applying the
general approach 
from this paper or
\cite{TSOS},  
the isotropy gro\-ups un\-der orthogonal similarity on skew-symmet\-ric or orthogonal matrices
are describ\-ed by 
equations involving a significant difference in comparison to 
(\ref{eqFYFIY}).
However, this problem is expected to be addressed in a future study.
\end{enumerate}
\end{remark}

%
%
%
%




\end{document}